\documentclass[10.5pt]{article}
\usepackage[leqno]{amsmath}
\usepackage{amsfonts}
\usepackage{graphicx}

\def\sqr#1#2{{\vcenter{\vbox{\hrule height.#2pt
              \hbox{\vrule width.#2pt height#1pt \kern#1pt \vrule width.#2pt}
              \hrule height.#2pt}}}}
\def\signed #1{{\unskip\nobreak\hfil\penalty50
              \hskip2em\hbox{}\nobreak\hfil#1
              \parfillskip=0pt \finalhyphendemerits=0 \par}}
\def\endpf{\signed {$\sqr69$}}

\def\3n{\negthinspace \negthinspace \negthinspace }
\def\2n{\negthinspace \negthinspace }
\def\1n{\negthinspace }

\def\dbE{\mathbb{E}}
\def\dbF{\mathbb{F}}

\def\dbP{\mathbb{P}}

\def\dbR{\mathbb{R}}

\def\={\buildrel \triangle \over =}

\def\ds{\displaystyle}

\def\ns{\noalign{\ss}}

\def\a{\alpha}
\def\b{\beta}

\def\d{\delta}
\def\e{\varepsilon}
\def\z{\zeta}

\def\l{\lambda}
\def\m{\mu}
\def\n{\nu}
\def\si{\sigma}
\def\t{\tau}
\def\f{\varphi}
\def\th{\theta}

\def\D{\Delta}
\def\Th{\Theta}

\def\O{\Omega}

\def\cF{{\cal F}}

\def\cH{{\cal H}}

\def\cM{{\cal M}}

\def\cU{{\cal U}}

\def\cl{{\cal l}}
\def\ss{\smallskip}
\def\ms{\medskip}

\def\q{\quad}
\def\qq{\qquad}
\def\hb{\hbox}

\def\lan{\mathop{\langle}}
\def\ran{\mathop{\rangle}}

\def\wt{\widetilde}

\def\cd{\cdot}

\def\ae{\hbox{\rm a.e.{ }}}
\def\as{\hbox{\rm a.s.{ }}}

\def\cl{\overline}

\def\deq{\mathop{\buildrel\D\over=}}

\def\({\Big (}
\def\){\Big )}
\def\[{\Big[}
\def\]{\Big]}
\def\bde{\begin{definition}}
\def\ede{\end{definition}}
\def\be{\begin{equation}}
\def\bel{\begin{equation}\label}
\def\ee{\end{equation}}
\def\bt{\begin{theorem}}
\def\et{\end{theorem}}
\def\bc{\begin{corollary}}
\def\ec{\end{corollary}}
\def\bl{\begin{lemma}}
\def\el{\end{lemma}}
\def\bp{\begin{proposition}}
\def\ep{\end{proposition}}
\def\bas{\begin{assumption}}
\def\eas{\end{assumption}}
\def\br{\begin{remark}}
\def\er{\end{remark}}
\def\ba{\begin{array}}
\def\ea{\end{array}}
\def\ed{\end{document}}

\def\square#1{\vbox{\hrule\hbox{\vrule height#1%
     \kern#1\vrule}\hrule}}
\def\rectangle#1#2{\vbox{\hrule\hbox{\vrule height#1%
     \kern#2\vrule}\hrule}}

\font\tenbb=msbm10 \font\sevenbb=msbm7 \font\fivebb=msbm5

\newfam\bbfam
\scriptscriptfont\bbfam=\fivebb \textfont\bbfam=\tenbb
\scriptfont\bbfam=\sevenbb

\newtheorem{lemma}{Lemma}[section]
\newtheorem{remark}{Remark}[section]

\newtheorem{theorem}{Theorem}[section]
\newtheorem{corollary}{Corollary}[section]

\newtheorem{definition}{Definition}[section]
\newtheorem{proposition}{Proposition}[section]
\newtheorem{assumption}{Assumption}[section]

\makeatletter
   
   \@addtoreset{equation}{section}
\makeatother

\title{\bf Optimal Control Problems of Forward-Backward\\
Stochastic Volterra Integral Equations}

\date{}

\author{
Yufeng Shi\footnote{ Institute for Financial Studies, School of
Mathematics, Shandong University, Jinan 250100, China. The research
of this author is partially supported by NNSF of China (Grants 11371226,
11071145 and 11231005), Foundation for Innovative Research Groups of
NNSF of China (Grant 11221061) and the 111 Project (Grant B12023).
},~~Tianxiao Wang\footnote{School of Mathematics, Sichuan
University, Chengdu 610065, China. The research of this author is
supported by National Basic Research Program of China (973 Program)
under grant 2011CB808002, the NNSF of China under grant 11231007 and
11301298, the PCSIRT under grant IRT1273 from Chinese Education
Ministry} ~~and~~Jiongmin Yong\footnote{Department of Mathematics,
University of Central Florida, Orlando, FL 32816, USA. The research
of this author is partially supported by NSF grant DMS-1007514. }}

\begin{document}

\maketitle

\centerline{(Dedicated to Professor Xunjing Li on His 80th Birth
Anniversary)}

\begin{abstract}
Optimal control problems of forward-backward stochastic Volterra
integral equations (FBSVIEs in short) are formulated and studied. A
general duality principle is established for linear backward
stochastic integral equation and linear stochastic Fredholm-Volterra
integral equation with mean-field. With the help of such a duality
principle, together with some other new delicate and subtle skills,
Pontryagin type maximum principles are proved for two optimal
control problems of FBSVIEs.

\end{abstract}

\bf Keywords: \rm Forward-backward stochastic Volterra integral
equations, adapted M-solution, stochastic maximum principle, duality
principle, stochastic Fredholm-Volterra integral equations

\ms

\bf AMS Mathematics subject classification. \rm 60H20, 93E20.

%% \linenumbers

%% main text

\section{Introduction}
%\label{sec:intro}

Let $(\O,\cF,\dbF,\dbP)$ be a complete filtered probability space on
which a standard one-dimensional Brownian motion $\{W(t),t\ge0\}$ is
defined with $\dbF=\{\cF_t\bigm|t\ge0\}$ being its natural
filtration augmented by all the $\dbP$-null sets. We point out that
assuming $W(\cd)$ to be one-dimensional is just for the simplicity
of presentation; Our results remain for the case of
multi-dimensional Brownian motions.

\ms

Let us start with a classical stochastic optimal control problem. To
this end, we consider the following controlled stochastic
differential equation (SDE, for short):
\bel{SDE}\left\{\2n\ba{ll}
\ns\ds dX(s)=b(s,X(s),u(s))ds+\si(s,X(s),u(s))dW(s),\qq s\in[0,T],\\
\ns\ds X(0)=x,\ea\right.\ee
with cost functional
\bel{cost0}J^0(x;u(\cd))=\dbE\Big[h(X(T))+\int_0^Tg(s,X(s),u(s))ds\Big],\ee
where $b,\si,h,g$ are some suitable maps, the {\it control} $u(\cd)$
is taken from some suitable set $\cU$, and the {\it state process}
$X(\cd)$ is valued in $\dbR^n$. In the above, components of $X(\cd)$
could be wealth, owned commodities/assets, inventory of products,
and some economic factors (interest rates, unemployment rate), and
so on. Also, $u(\cd)$ can be regarded as some kind of investment,
production effort, labor force, transaction of assets, etc. Then our
classical optimal control can be stated as follows.

\ms

\bf Problem (C)$^0$. \rm For given $x\in\dbR^n$, find a $\bar
u(\cd)\in\cU$, called an {\it optimal control}, such that
\bel{}J^0(x;\bar u(\cd))=\inf_{u(\cd)\in\cU}J(x;u(\cd)).\ee

\ms

Standard results for the above Problem (C)$^0$ can be found, say, in
\cite{Yong-Zhou 1999}.

\ms

Now let us make some further analysis on the equation (\ref{SDE}) which
can be written as follows:
\bel{1.4}X(t)=x+\int_0^tb(s,X(s),u(s))ds+\int_0^t\si(s,X(s),u(s))dW(s),\qq
t\in[0,T].\ee
Although the state equation has the above integral form, it is still
{\it memoryless}, in the sense that the increment $X(t+\D t)-X(t)$
of the state on $[t,t+\D t]$ only depends on the local ``driving
force'' $\big\{\big(b(s,X(s),u(s)),\si(s,X(s)$,
$u(s))\big),~s\in[t,t+\D t]\big\}$:
$$X(t+\D t)-X(t)=\int_t^{t+\D t}b(s,X(s),u(s))ds+\int_t^{t+\D t}\si(s,X(s),u(s))dW(s).$$
Whereas, in reality, memory or long-term dependence often exists. In
another word, the increment $X(t+\D t)-X(t)$ of the state on
$[t,t+\D t]$ might depend on the ``driving force'' of
non-infinitesimal time duration, say, $[t-\t,t]$, for some $\t>0$.
For example, the current production level usually depends on some
renovation of production equipments some time ago, the profit of
investment usually depends on the transactions some time before, the
air pollution is caused by some bad production strategies some years
ago, etc. To model various possible situations with memory, instead
of (\ref{1.4}), we may consider the following controlled (forward)
stochastic Volterra integral equation (FSVIE, for short):
\bel{FSVIE1}X(t)=\f(t)+\int_0^tb(t,s,X(s),u(s))ds+\int_0^t\si(t,s,X(s),
u(s))dW(s),\q t\in[0,T].\ee
Unlike (\ref{1.4}), due to the dependence of
$\big(b(t,s,x,u),\si(t,s,x,u)\big)$ on $t$, even for the case
$\f(t)\equiv x$, we have
$$\ba{ll}
\ns\ds X(t+\D t)-X(t)=\int_t^{t+\D t}b(t+\D
t,s,X(s),u(s))ds+\int_t^{t+\D
t}\si(t+\D t,s,X(s),u(s)dW(s)\\
\ns\ds\qq\qq\qq\qq\qq+\int_0^t\big[b(t+\D
t,s,X(s),u(s))-b(t,s,X(s),u(s))\big]ds\\
\ns\ds\qq\qq\qq\qq\qq+\int_0^t\big[\si(t+\D
t,s,X(s),u(s))-\si(t,s,X(s),u(s))\big]dW(s),\ea$$
which depends not only on the values of the ``driving force'' in
$[t,t+\D t]$, but also on those in the whole interval $[0,t]$ up to
the current time $t$. Therefore, with suitable choices of $b$ and
$\si$, it is possible to model certain memory effects through
(\ref{FSVIE1}). Based on these arguments, one can use deterministic
or stochastic Volterra integral equations to describe some economic
models, see \cite{Duffie-Huang 1986,Pardoux-Protter
1990,Kamien-Muller 1976}, for examples.

\ms

On the other hand, let us turn to the cost functional (\ref{cost0}). As
we know, stochastic differential utility (SDU, for short) introduced
by Duffie--Epstein (\cite{Duffie-Epstein 1992}) can be represented
by {\it backward stochastic differential equations} (BSDEs, for
short). More precisely, if $C(\cd)$ is a consumption process and
$\xi$ is a payoff at the terminal time $T$, then an SDU process
$Y(\cd)$ for the pair $(C(\cd),\xi)$ can be modeled by the
following:
\bel{SDU}Y(t)=\dbE\Big[\xi+\int_t^Tg(s,Y(s),C(s))ds\Bigm|\cF_t\Big],\qq
t\in[0,T],\ee
for some suitable map $g$. This can also be regarded as a {\it
dynamic risk measure process} associated with the pair
$(\xi,C(\cd))$. It turns out that (\ref{SDU}) admits the following
equivalent form,
\bel{BSDE1}\ba{ll}
\ns\ds Y(t)=\xi+\int_t^Tg(s,Y(s),C(s))ds-\int_t^TZ(s)dW(s),\qq
s\in[0,T], \ea\ee
which is a BSDE whose solution is a pair $(Y(\cd),Z(\cd))$ of
$\dbF$-adapted processes (\cite{Pardoux-Protter 1990, Ma-Yong 1999,
Yong-Zhou 1999}). Note that if we let
$$\xi=h(X(T)),\qq u(\cd)=C(\cd),$$
then the cost functional (\ref{cost0}) admits the following
representation:
$$J^0(x;u(\cd))=Y(0).$$
 Thanks to the further development of BSDEs,
one could extend the SDU theory via more general BSDEs, namely, one
may define an SDU process $(Y(\cd),Z(\cd))$ as the adapted solution
to the following general BSDE:
\bel{BSDE2}\left\{\2n\ba{ll}
\ns\ds dY(t)=-g(t,Y(t),Z(t),C(t))dt+Z(t)dW(t),\qq t\in[0,T],\\
\ns\ds Y(T)=\xi,\ea\right.\ee
whose equivalent integral form reads
$$Y(t)=\xi+\int_t^Tg(s,Y(s),Z(s),C(s))ds-\int_t^TZ(s)dW(s),\qq
t\in[0,T].$$
Although it is very general, the above is still challenged by the
following two aspects: (i) The terminal payoff/cost $\xi$ is
time-independent; (ii) the ``running utility/cost'' is of memoryless
feature. These two lead to the {\it time-consistency} of the utility
process $Y(\cd)$, which is a little too ideal. In reality,
substantial evidence (see \cite{Ainslie 1992}, for example) shows
that people in the real life are more concerned (or impatient) about
the choices (or decisions) for the immediate future, but are more
rational (or patient) when facing long-term alternatives. Such a
phenomenon is just one particular case of time inconsistency.
Therefore inspired by the theory of backward stochastic Volterra
integral equations (BSVIEs, for short) (\cite{Lin 2002, Wang-Shi
2010, Wang-Yong 2013, Wang-Zhang 2007, Yong 2006, Yong 2007, Yong
2008}), we could modify the above into the following form, taking
into account of our controlled FSVIE:
\bel{1.9}\ba{ll}
\ns\ds
Y(t)=\psi(t,X(t),X(T))+\int_t^Tg(t,s,X(t),X(s),Y(s),Z(t,s),Z(s,t),u(s))ds\\
\ns\ds\qq\qq\qq-\int_t^TZ(t,s)dW(s),\qq t\in[0,T].\ea\ee
where $(Y(\cd),Z(\cd\,,\cd))$ is a so-called {\it adapted
M-solution} of the above (see \cite{Yong 2008}). In the above, we
see that both $X(t)$ and $X(T)$ appear in the free-term
$\psi(t,X(t),X(T))$. A motivation of that is the following: Suppose
$X(\cd)$ represent the production level process of certain product.
One expects that the terminal level should be within a certain range
determined by the current level, due to the limitation of resource,
manpower, machine capacity, market demand/price, etc. Some similar
explanations can be made for the appearance of both $X(t)$ and
$X(s)$ in the integrand.

\ms

Motivated by the above arguments, in this paper we study the
following controlled forward-backward stochastic Volterra integral equations
(FBSVIEs, for short):
\bel{equation}\left\{\2n\ba{ll}
\ns\ds
X(t)=\f(t)+\int_0^tb(t,s,X(s),u(s))ds+\int_0^t\si(t,s,X(s),u(s))
dW(s), \\
\ns\ds
Y(t)=\psi(t,X(t),X(T))+\int_t^Tg(t,s,X(t),X(s),Y(s),Z(t,s),Z(s,t),u(s))ds\\
\ns\ds\qq\qq\qq-\int_t^TZ(t,s)dW(s),\qq t\in[0,T].\ea\right.\ee
We call $(X(\cd),Y(\cd),Z(\cd\,,\cd))$ the {\it state} and $u(\cd)$
the {\it control}. In such a system, $X(\cd)$ and $Y(\cd)$ can be
regarded as the {\it portfolio} process and the {\it dynamic risk}
process, respectively. To introduce the cost functional, we need to
separate two cases.

\ms

First of all, if the {\it generator} $g(\cd)$ of the BSVIE in
(\ref{equation}) is independent of $Z(s,t)$, then the state equation
reads:
\bel{equation2}\left\{\2n\ba{ll}
\ns\ds
X(t)=\f(t)+\int_0^tb(t,s,X(s),u(s))ds+\int_0^t\si(t,s,X(s),u(s))
dW(s), \\
\ns\ds
Y(t)=\psi(t,X(t),X(T))+\int_t^Tg(t,s,X(t),X(s),Y(s),Z(t,s),u(s))ds\\
\ns\ds\qq\qq\qq-\int_t^TZ(t,s)dW(s),\qq t\in[0,T].\ea\right.\ee
In this case, under some mild conditions, for any control $u(\cd)$,
there exists a unique triplet $(X(\cd),Y(\cd)$, $Z(\cd\,,\cd))$,
called the {\it adapted solution} to (\ref{equation2}), such that
\bel{}\left\{\2n\ba{ll}
\ns\ds s\mapsto(X(s),Y(s))\hb{ is $\dbF$-adapted on $[0,T]$},\\
\ns\ds s\mapsto Z(t,s)\hb{ is $\dbF$-adapted on $[t,T]$,}\q\forall
t\in[0,T),\ea\right.\ee
and (\ref{equation2}) is satisfied in the usual It\^o's sense.
Moreover, $Y(\cd)\in C_\dbF(0,T;L^2(\O;\dbR^n))$. Therefore, $Y(0)$
is well-defined and for such a case, we may introduce the following
cost functional:
\bel{cost1}J_1(u(\cd))=\dbE\[h(X(T),Y(0))
+\int_0^T\int_t^Tf(t,s,X(s),Y(s),Z(t,s),u(s))dsdt\].\ee
Note that in the above case, the process $Z(t,s)$ is only defined for
$0\le t\le s\le T$.

\ms

On the other hand, if the generator $g(\cd)$ depends on $Z(s,t)$,
then, by \cite{Yong 2008}, under some suitable conditions, there
exists a unique triplet $(X(\cd),Y(\cd),Z(\cd\,,\cd))$, called the
{\it adapted M-solution} to (\ref{equation}), such that
\bel{}\left\{\2n\ba{ll}
\ns\ds s\mapsto(X(s),Y(s))\hb{ is $\dbF$-adapted on $[0,T]$},\\
\ns\ds s\mapsto Z(t,s)\hb{ is $\dbF$-adapted on $[0,T]$,}\q\ae
t\in[0,T],\ea\right.\ee
and in addition to (\ref{equation}) being satisfied in the usual
It\^o's sense, one also has
\bel{1.15}Y(t)=\dbE Y(t)+\int_0^tZ(t,s)dW(s),\qq\as\2n,~\ae
t\in[0,T].\ee
Different from the first case, in this second case, the process $Z(t,s)$
is defined on $[0,T]^2$, and the additional relation (\ref{1.15})
holds. In this second case, due to the fact that $t\mapsto Z(s,t)$
(for $t\le s$) is not necessarily continuous, we could not expect
the continuity of $t\mapsto Y(t)$. Therefore, $Y(0)$ might not be
well-defined in general. Consequently, the corresponding cost
functional should not contain the term like $h(X(T),Y(0))$ as in
$J_1(u(\cd))$. Fortunately, a comparison theorem found in
\cite{Wang-Yong 2013} suggests that in the current case, it might be
more proper to use $\dbE\int_0^TY(s)ds$ as an alternative for
$Y(0)$. Hence, we propose the following cost functional:
\bel{cost2}J_2(u(\cd))=\dbE\[h\(X(T),\dbE\int_0^TY(s)ds\)
+\int_0^T\int_0^Tf(t,s,X(s),Y(s),Z(t,s),u(s))dsdt\].\ee

\ms

From the above, we see that one can formulate two different optimal
control problems. In this paper, we will establish Pontryagin type
maximum principles for the optimal control problems corresponding to
the above two settings. It is known that in deriving maximum
principle, besides the suitable variation of the state equation and
cost functional, the key is to have a duality principle. The major
contribution of this paper is the discovery of a duality principle for
general linear BSVIEs, which is a significant extension of that
presented in \cite{Yong 2008}. It turns out that our new duality
principle involves a special type of stochastic Fredholm-Volterra
integral equations, and we are able to obtain its solvability under
natural conditions. It is worthy of pointing out that in contract
with SDE case (\cite{Peng 1990, Peng 1993, Yong-Zhou 1999}), we need
to carry out all the calculations without differentiation due to the
lack of It\^o's formula for stochastic integral equations.

\ms

The rest of this paper is organized as follows. In Section 2, some
basic results concerning BSVIEs are recalled. In Section 3, we state
two maximum principles for controlled FBSVIEs. A general dual
principle for linear BSVIEs is established in Section 4. Then the
stated maximum principles are proved in Section 5. Section 6
concludes the paper.

\ms

\section{Results for BSVIEs Revisited}

In this section, we are going to recall some relevant results for
BSVIEs. To this end, let us first introduce some spaces. For
$H=\dbR^n,\dbR^{n\times m}$, etc., we denote its norm by $|\cd|$.
For $0\le s<t\le T$, define
$$\ba{ll}
\ns\ds L^2_{\cF_t}(\O;H)=\Big\{\xi:\O\to H\bigm|\xi\hb{ is
$\cF_t$-measurable, }\dbE|\xi|^2<\infty\Big\},\\
\ns\ds L^2_{\cF_T}(s,t;H)=\Big\{X:[s,t]\times\O\to H\bigm|X(\cd)\hb{
is $\cF_T$-measurable},~\dbE\int_s^t|X(r)|^2dr<\infty\Big\},\\
\ns\ds L^2_{\cF_T}\big(\O;C(s,t;H)\big)=\Big\{X:[s,t]\times\O\to
H\bigm|X(\cd)\hb{ is $\cF_T$-measurable, has continuous paths,
}\\
\ns\ds\qq\qq\qq\qq\qq\qq\qq\qq\qq\qq\qq\dbE\(\sup_{r\in[s,t]}|X(r)|^2\)<\infty\Big\},\\
\ns\ds C_{\cF_T}\big(s,t;L^2(\O;H)\big)\1n=\1n\Big\{X\1n:\1n[s,t]\to
L^2_{\cF_T}(\O;H)\bigm|X(\cd)\hb{ is
continuous,}\sup_{r\in[s,t]}\dbE|X(r)|^2<\infty\Big\}.\ea$$
Also, we define
$$L^2_\dbF(s,t;H)=\Big\{X(\cd)\in L^2_{\cF_T}(s,t;H)\bigm|X(\cd)\hb{ is $\dbF$-adapted}\Big\}.$$
The spaces $L^2_\dbF\big(\O;C(s,t;H)\big)$ and
$C_\dbF\big(s,t;L^2(\O;H)\big)$ can be defined in the same way. It
should be pointed out that
$$\ba{ll}
\ns\ds L^2_{\cF_T}(\O;C(s,t;H))\subseteq
C_{\cF_T}(s,t;L^2(\O;H)),\\
\ns\ds L^2_\dbF(\O;C(s,t;H))\subseteq C_\dbF(s,t;L^2(\O;H)),\ea$$
and the equalities do not hold in general. Further, we denote
$$\D=\Big\{(t,s)\in[0,T]^2\Bigm|t\le s\Big\},\q\D^*=\Big\{(t,s)\in[0,T]^2\Bigm|t\ge s\Big\}\equiv
\cl{\D^c},$$
and let
$$\ba{ll}
\ns\ds L^2_\dbF(\D;H)=\Big\{Z:\D\times\O\to H\bigm|s\mapsto
Z(t,s)\hb{ is $\dbF$-adapted on $[t,T]$, $\ae t\in[0,T]$,}\\
\ns\ds\qq\qq\qq\qq\qq\qq\qq\|Z(\cd\,,\cd)\|_{L^2_\dbF(\D;H)}^2\equiv\dbE\int_0^T
\int_t^T|Z(t,s)|^2dsdt<\infty\Big\},\ea$$
$$\ba{ll}
\ns\ds L^2(0,T;L^2_\dbF(0,T;H))=\Big\{Z:[0,T]^2\times\O\to
H\bigm|s\mapsto Z(t,s)\hb{ is $\dbF$-adapted}\hb{ on $[0,T]$, $\ae
t\in[0,T]$,}\\
\ns\ds\qq\qq\qq\qq\qq\qq\qq\qq\|Z(\cd\,,\cd)\|_{L^2_\dbF(0,T;L^2(0,T;H))}^2\equiv
\dbE\int_0^T\int_0^T|Z(t,s)|^2dsdt<\infty\Big\}.\ea$$
We denote
$$\left\{\2n\ba{ll}
\ns\ds\cH^2_\D[0,T]=L^2_\dbF(0,T;\dbR^m)\times L^2_\dbF(\D;\dbR^m),\\
\ns\ds\cH^2[0,T]=L^2_\dbF(0,T;\dbR^m)\times
L^2\big(0,T;L^2_\dbF(0,T;\dbR^m)\big).\ea\right.$$
Further, we let $\cM^2[0,T]$ be the set of all
$(y(\cd),z(\cd\,,\cd))\in\cH^2[0,T]$ such that
$$y(t)=\dbE y(t)+\int_0^tz(t,s)dW(s),\q\as\2n,~\ae t\in[0,T].$$
Clearly, $\cM^2[0,T]$ is a closed subspace of $\cH^2[0,T]$. Also,
for any $(y(\cd),z(\cd\,,\cd))\in\cM^2[0,T]$, we have
\bel{2.1}\dbE|y(t)|^2=\big|\dbE
y(t)\big|^2+\dbE\int_0^t|z(t,s)|^2ds\ge\dbE\int_0^t|z(t,s)|^2ds,\qq\ae
t\in[0,T].\ee
The above implies that for any $\b\ge0$, there exists a constant
$K>0$ depending on $\b$ such that
\bel{2.2}\ba{ll}
\ns\ds\|(y(\cd),z(\cd\,,\cd))\|_{\cH^2[0,T]}^2\equiv\dbE\[\int_0^T|y(t)|^2dt
+\int_0^T\int_0^T|z(t,s)|^2dsdt\]\\
\ns\ds\le2\dbE\[\int_0^T|y(t)|^2dt
+\int_0^T\int_t^T|z(t,s)|^2dsdt\]\equiv2\|(y(\cd),z(\cd\,,\cd))\|_{
\cH^2_\D[0,T]}^2\\
\ns\ds\le K\dbE\[\int_0^Te^{\b t}|y(t)|^2dt+\int_0^Te^{\b
t}\int_t^T|z(t,s)|^2dsdt\]
\equiv K\|(y(\cd),z(\cd\,,\cd))\|_{\cM^2_\b[0,T]}^2\\
\ns\ds\le K\|(y(\cd),z(\cd\,,\cd))\|_{\cH^2[0,T]}^2,\ea\ee
which means that $\|\,\cd\,\|_{\cM^2_\b[0,T]}$ is an equivalent norm
of $\|\,\cd\,\|_{\cH^2[0,T]}$ on $\cM^2[0,T]$. Next, we let
$\wt\cM^2[0,T]$ be the set of all
$(y(\cd),z(\cd\,,\cd))\in\cM^2[0,T]$ such that $t\mapsto z(s,t)$ is
continuous on $[0,s]$. Following the ideas in, for example, Lemma
2.4 (p.133) and Proposition 2.6 (p.134) of \cite{Karatzas-Shreve
1988}, one can see that any measurable process in $\cM^2[0,T]$ can
be approximated by the one in $\wt\cM^2[0,T]$. Therefore
$\wt\cM^2[0,T]$ is dense in $\cM^2[0,T]$.

 \ms

\rm

We now consider the following two types of BSVIE:
\bel{2.3}Y(t)=\psi(t)+\int_t^Tg(t,s,Y(s),Z(t,s),Z(s,t))ds-\int_t^TZ(t,s)dW(s),\q\as\2n,~\ae
t\in[0,T].\ee
\bel{2.4}Y(t)=\psi(t)+\int_t^Tg(t,s,Y(s),Z(t,s))ds-\int_t^TZ(t,s)dW(s),\q\as\2n,~\ae
t\in[0,T],\ee
Note that (\ref{2.4}) is a special case of (\ref{2.3}) with the
generator $g(\cd)$ independent of $Z(s,t)$. We recall the following
definition.

\ms

\bf Definition 2.2. \rm (i) A pair
$(Y(\cd),Z(\cd\,,\cd))\in\cH^2_\D[0,T]$ is called an {\it adapted
solution} to BSVIE (\ref{2.4}) if for almost every $t\in[0,T]$,
$Y(t)$ is $\cF_t$-measurable, $s\mapsto Z(t,s)$ is $\dbF$-adapted on
$[t,T]$, and (\ref{2.4}) is satisfied in the usual It\^o's sense.

\ms

(ii) A pair $(Y(\cd),Z(\cd\,,\cd))\in\cM^2[0,T]$ is called an {\it
adapted M-solution} to BSVIE (\ref{2.3}) if for almost every
$t\in[0,T]$, $Y(t)$ is $\cF_t$-measurable, $s\mapsto Z(t,s)$ is
$\dbF$-adapted on $[0,T]$, and (\ref{2.3}) is satisfied in the usual
It\^o's sense.

\ms

We point out that for BSVIE (\ref{2.4}) the values $Z(t,s)$ of
$Z(\cd\,,\cd)$ with $t\ge s$ are irrelevant. Hence, adapted
solutions $(y(\cd),z(\cd\,,\cd))$ of (\ref{2.4}) need only belong to
$\cH^2_\D[0,T]$.

\ms

The following is a hypothesis for the coefficients of BSVIE
(\ref{2.3}).

\ms

{\bf(H1)} Let
$g:[0,T]^2\times\dbR^m\times\dbR^m\times\dbR^m\times\O\to\dbR^m$ be
measurable such that $(y,z,z')\mapsto g(t,s,y,z,z')$ is uniformly
Lipschitz, and
\bel{}\dbE\int_0^T\left(\int_0^T|g(t,s,0,0,0)|ds\right)^2dt\le
\infty,\quad~y,z,z'\in\dbR^m,~\ee

\ms

For BSVIE (\ref{2.4}), we introduce the following hypothesis:

\ms

{\bf(H2)} Let $g:[0,T]^2\times\dbR^m\times\dbR^m\times\O\to\dbR^n$
be measurable such that $(y,z)\mapsto g(t,s,y,z)$ is uniformly
Lipschitz,
\bel{}\dbE\int_0^T\left(\int_0^T|g(t,s,0,0)|ds\right)^2dt\le
\infty,\ee
and there exists a modulus of
continuity $\rho:[0,\infty)\to[0,\infty)$ (a continuous and monotone
increasing function with $\rho(0)=0$) such that
\bel{}|g(t,s,y,z)-g(t',s,y,z)|\le\rho(|t-t'|)(|y|+|z|),\qq\forall
t,t',s\in[0,T],~y,z\in\dbR^m.\ee

\ms

We have the following result concerning BSVIEs
(\ref{2.3})--(\ref{2.4}). One can essentially find proofs of such a
result in \cite{Yong 2006,Yong 2008,Shi-Wang 2012,Shi-Wang-Yong
2013}. For convenience, we provide a direct proof.

\ms

\bf Theorem 2.3. \sl {\rm(i)} Let {\rm(H1)} hold. Then for any
$\psi(\cd)\in L^2_{\cF_T}(0,T;\dbR^m)$, BSVIE $(\ref{2.3})$ admits a
unique adapted M-solution $(Y(\cd),Z(\cd\,,\cd))\in\cH^2[0,T]$, and
the following holds:
\bel{2.7}\ba{ll}
\ns\ds\dbE\int_0^T|Y(t)|^2ds+\dbE\int_0^T\int_0^T|Z(t,s)|^2dsdt\\
\ns\ds\le
K\[\dbE\int_0^T|\psi(t)|^2dt+\int_0^T\(\int_0^T|g(t,s,0,0,0)|ds\)^2dt\].\ea\ee
If $\psi_i(\cd)\in L^2_{\cF_T}(0,T;\dbR^m)$, $g_i(\cd)$ satisfies
(H1) and $(Y_i(\cd),Z_i(\cd\,,\cd))$ are corresponding adapted
M-solutions, then
\bel{stability}\ba{ll}
\ns\ds\dbE\int_0^T|Y_1(t)-Y_2(t)|^2dt+\dbE\int_0^T\int_0^T|Z_1(t,s)-Z_2(t,s)|^2dsdt\\
\ns\ds\le
K\[\int_0^T|\psi_1(t)-\psi_2(t)|^2dt\\
\ns\ds\qq+\int_0^T\(\int_0^T|g_1(t,s,Y_1(s),Z_1(t,s),Z_1(s,t))
-g_2(t,s,Y_1(s),Z_1(t,s),Z_1(s,t))|ds\)^2dt\].\ea\ee

\ms

{\rm(ii)} Let {\rm(H2)} hold. Then for any $\psi(\cd)\in
C_{\cF_T}(0,T;L^2(\O;\dbR^m))$, BSVIE $(\ref{2.4})$ admits a unique
adapted solution $(Y(\cd),Z(\cd,\cd))\in
C_\dbF(0,T;L^2(\O;\dbR^m))\times L^2_{\dbF}(\D;\dbR^m)$ such that
\bel{2.9}\ba{ll}
\ns\ds\sup_{t\in[0,T]}\dbE|Y(t)|^2+\sup_{t\in[0,T]}\dbE\int_t^T|Z(t,s)|^2ds\\
\ns\ds\le K\Big\{\sup_{t\in[0,T]}\dbE|\psi(t)|^2
+\sup_{t\in[0,T]}\dbE\Big(\int_t^T|g(t,s,0,0)|ds\Big)^2\Big\}.\ea\ee
If $g_i(\cd)$ satisfies {\rm(H2)}, $\psi_i(\cd)\in
C_{\cF_T}(0,T;L^2(\Omega;\dbR^m))$, and $(Y_i(\cd),Z_i(\cd,\cd))\in
C_\dbF(0,T;L^2(\O;\dbR^m))\times L^2_{\dbF}(\D;\dbR^m)$ is the
unique adapted solution of BSVIE $(\ref{2.4})$ corresponding to
$\(g_i(\cd),\psi_i(\cd)\),$ then
\bel{2.10}\ba{ll}
\ns\ds\sup_{t\in[0,T]}\dbE|Y_1(t)-Y_2(t)|^2+\sup_{t\in[0,T]}\dbE\int_t^T|Z_1(t,s)-Z_2(t,s)|^2ds\\
\ns\ds\le
K\Big\{\sup_{t\in[0,T]}\dbE|\psi_1(t)-\psi_2(t)|^2\\
\ns\ds\qq+\sup_{t\in[0,T]}\dbE\Big(\int_t^T|g_1(t,s,Y_1(s),Z_1(t,s))-
g_2(t,s,Y_1(s),Z_1(t,s))|ds\Big)^2\Big\}.\ea\ee

\rm

\ms

\it Proof. \rm First, we let $\psi(\cd)\in
C_{\cF_T}([0,T];L^2(\O;\dbR^m))$, and assume that $t\mapsto
g(t,s,y,z,z')$ is continuous. For any
$(y(\cd),z(\cd\,,\cd))\in\wt\cM^2[0,T]$, consider the following BSDE
parameterized by $t\in[0,T]$:
\bel{BSDE}\eta(t,r)=\psi(t)+\int_r^Tg(t,s,y(s),\z(t,s),z(s,t))ds-\int_r^T\z(t,s)dW(s),\q
r\in[0,T].\ee
Under our conditions, the above BSDE admits a unique adapted
solution
$$(\eta(t,\cd),\z(t,\cd))\in L^2_\dbF(\O;C([0,T];\dbR^m))
\times L^2_\dbF(\D;\dbR^m),$$
and the following holds:
\bel{2.11}\ba{ll}
\ns\ds\dbE\[\sup_{r\in[t,T]}|\eta(t,r)|^2+\int_t^T|\z(t,s)|^2ds\]\le
K\dbE\[|\psi(t)|^2+\(\int_t^T|g(t,s,y(s),0,z(s,t))|ds\)^2\]\\
\ns\ds\le
K\dbE\Big\{|\psi(t)|^2+\int_t^T|y(s)|^2ds+\(\int_t^T|z(s,t)|ds\)^2+
\(\int_t^T|g(t,s,0,0,0)|ds\)^2\Big\}.\ea\ee
We now let
\bel{YZ1}Y(t)=\eta(t,t),\q t\in[0,T],\qq
Z(t,s)=\z(t,s),\qq(t,s)\in\D,\ee
and define $Z(t,s)$ with $t\ge s$ through the following:
\bel{YZ2}Y(t)=\dbE Y(t)+\int_0^tZ(t,s)dW(s),\qq t\in[0,T].\ee
Then $(Y(\cd),Z(\cd\,,\cd))$ is an adapted M-solution to the
following BSVIE:
$$Y(t)=\psi(t)+\int_t^Tg(t,s,y(s),Z(t,s),z(s,t))ds-\int_t^TZ(t,s)dW(s),\q t\in[0,T],$$
and from (\ref{2.11}), one has
\bel{2.15}\ba{ll}
\ns\ds\dbE\[|Y(t)|^2+\int_t^T|Z(t,s)|^2ds\]\\
\ns\ds\le
K\dbE\Big\{|\psi(t)|^2+\(\int_t^T|g_0(t,s)|ds\)^2+\int_t^T|y(s)|^2ds+\(\int_t^T|z(s,t)|ds\)^2\Big\},\ea\ee
where
$$g_0(t,s)=g(t,s,0,0,0).$$
Consequently, noting (\ref{2.1}),
$$\ba{ll}
\ns\ds\|(Y(\cd),Z(\cd\,,\cd))\|^2_{\cH^2[0,T]}\equiv\dbE\Big\{\int_0^T|Y(t)|^2dt
+\int_0^T\int_0^T|Z(t,s)|^2dsdt\Big\}\\
\ns\ds\le\dbE\Big\{2\int_0^T|Y(t)|^2dt
+\int_0^T\int_t^T|Z(t,s)|^2dsdt\Big\}\\
\ns\ds\le
K\dbE\Big\{\(\int_t^T\2n|g_0(t,s)|ds\)^2\3n+\1n\int_0^T\2n|\psi(t)|^2dt+\1n\int_0^T\2n\int_t^T|y(s)|^2dsdt
\1n+\1n\int_0^T\2n\(\int_t^T\2n|z(s,t)|ds\)^2dt\Big\}\\
\ns\ds\le
K\dbE\Big\{\(\int_t^T|g_0(t,s)|ds\)^2+\int_0^T|\psi(t)|^2dt+\int_0^T|y(t)|^2dt+\int_0^T\int_0^t|z(t,s)|^2dsdt\Big\}\\
\ns\ds\le
K\dbE\Big\{\(\int_t^T|g_0(t,s)|ds\)^2+\int_0^T|\psi(t)|^2dt+\int_0^T|y(t)|^2dt\Big\}\\
\ns\ds\le
K\Big\{\(\int_t^T|g_0(t,s)|ds\)^2+\|\psi(\cd)\|^2_{L^2_{\cF_T}(0,T;\dbR^m)}+\|(y(\cd),z(\cd\,,\cd))\|_{\cM^2[0,T]}\Big\}.\ea$$
Hence, if we define
$\Th(y(\cd),z(\cd\,,\cd))=(Y(\cd),Z(\cd\,,\cd))$, then $\Th$ maps
from $\wt\cM^2[0,T]$ to $\cM^2[0,T]$. Note that such a map also
depends on the choice of the {\it free term} $\psi(\cd)$ and the
generator $g(\cd)$. We now show that the mapping $\Th$ can be
extended to $\cM^2[0,T]$ and the extension, still denoted by $\Th$,
is contractive and has a stability property with respect to
$(\psi(\cd),g(\cd))$. To this end, we take any $\psi_i(\cd)\in
C_{\cF_T}([0,T];L^2(\O;\dbR^m))$, $g_i(\cd)$ satisfies (H1) with
$t\mapsto g_i(t,s,y,z,z')$ being continuous, and
$(y_i(\cd),z_i(\cd\,,\cd))\in\wt\cM^2[0,T]$ ($i=1,2$). Let
$(\eta_i(\cd),\z_i(\cd\,,\cd))$ be the adapted solution of BSDE
(\ref{BSDE}) with $(\psi(\cd),g(\cd))$ and $(y(\cd),z(\cd\,,\cd))$
replaced by $(\psi_i(\cd),g_i(\cd))$ and
$(y_i(\cd),z_i(\cd\,,\cd))$, respectively. By the stability of
adapted solutions to BSDEs, one has
$$\ba{ll}
\ns\ds\dbE\[\sup_{r\in[t,T]}|\eta_1(t,r)-\eta_2(t,r)|^2+\int_t^T|\z_1(t,s)-\z_2(t,s)|^2ds\]\\
\ns\ds\le
K\dbE\Big\{|\psi_1(t)-\psi_2(t)|^2+\(\int_t^T|g_1(t,s,y_1(s),\z_1(t,s),z_1(s,t))\\
\ns\ds\qq\qq\qq\qq\qq\qq\qq\qq-g_2(t,s,y_2(s),\z_1(t,s),z_2(s,t))|ds\)^2\Big\}\\
\ns\ds\le K\dbE\Big\{|\psi_1(t)-\psi_2(t)|^2+\(\int_t^T|\d g(t,s)|ds\)^2\\
\ns\ds\qq+\int_t^T|y_1(s)-y_2(s)|^2ds+\(\int_t^T|z_1(s,t)-z_2(s,t)|ds\)^2\Big\},
\ea$$
with
$$\d g(t,s)=g_1(t,s,y_1(s),\z_1(t,s),z_1(s,t))-g_2(t,s,y_1(s),\z_1(t,s),z_1(s,t)).$$
Then it by defining $(Y_i(\cd),Z_i(\cd\,,\cd))$ similar to
(\ref{YZ1})--(\ref{YZ2}), one has
\bel{2.16}\ba{ll}
\ns\ds\dbE\[|Y_1(t)-Y_2(t)|^2+\int_t^T|Z_1(t,s)-Z_2(t,s)|^2ds\]\\
\ns\ds\le K\dbE\[|\psi_1(t)-\psi_2(t)|^2+\(\int_t^T|\d
g(t,s)|ds\)^2\\
\ns\ds\qq+\int_t^T|y_1(s)-y_2(s)|^2ds+\(\int_t^T|z_1(s,t)-z_2(s,t)|ds\)^2\].\ea\ee
Making use of a relation similar to (\ref{2.1}) for
$(y_1(\cd)-y_2(\cd),z_1(\cd\,,\cd)-z_2(\cd\,,\cd))$, we have
\bel{2.17}\ba{ll}
\ns\ds\dbE\int_0^Te^{\b t}\[\int_t^T|y_1(s)-y_2(s)|^2ds+\(\int_t^T|z_1(s,t)-z_2(s,t)|ds\)^2\]dt\\
\ns\ds=\dbE\[\int_0^T\3n\(\1n\int_0^s\3n e^{\b t}dt\)|y_1(s)\1n-\1n
y_2(s)|^2ds\1n+\2n\int_0^T\3n e^{\b
t}\(\1n\int_t^T\3n e^{-{\b\over 2}s}e^{{\b\over 2}s}|z_1(s,t)\1n-\1n z_2(s,t)|ds\)^2dt\]\\
\ns\ds\le\dbE\[{1\over\b}\2n\int_0^T\3n e^{\b t}|y_1(t)\1n-\1n
y_2(t)|^2dt\1n+\3n\int_0^T\3n e^{\b t}\(\1n\int_t^T\3n e^{-\b
s}ds\)\(\int_t^T\3n e^{\b s}
|z_1(s,t)-z_2(s,t)|^2ds\)dt\]\\
\ns\ds\le{1\over\b}\dbE\int_0^Te^{\b
t}|y_1(t)-y_2(t)|^2dt+{1\over\b}\dbE\int_0^T\int_0^se^{\b
s}|z_1(s,t)-z_2(s,t)|^2dtds\\
\ns\ds\le{1\over\b}\dbE\int_0^Te^{\b
t}|y_1(t)-y_2(t)|^2dt+{1\over\b}\dbE\int_0^Te^{\b
s}\int_0^s|z_1(s,t)-z_2(s,t)|^2dtds\\
\ns\ds\le{2\over\b}\dbE\int_0^Te^{\b
t}|y_1(t)-y_2(t)|^2dt\le{2\over\b}\|(y_1(\cd),z_1(\cd\,,\cd))-(y_2(\cd),z_2(\cd\,,\cd))\|^2_{\cM^2_\b[0,T]}.\ea\ee
Consequently,
\bel{2.18}\ba{ll}
\ns\ds\|(Y_1(\cd),Z_1(\cd\,,\cd))-(Y_2(\cd),Z_2(\cd\,,\cd))\|^2_{\cM^2_\b[0,T]}\\
\ns\ds\equiv
\int_0^Te^{\b t}\dbE\[|Y_1(t)-Y_2(t)|^2+\int_t^T|Z_1(t,s)-Z_2(t,s)|^2ds\]dt\\
\ns\ds\le K\int_0^Te^{\b
t}\dbE\[|\psi_1(t)-\psi_2(t)|^2+\(\int_t^T|\d g(t,s)|ds\)^2\\
\ns\ds\qq+\int_t^T|y_1(s)-y_2(s)|^2ds+\(\int_t^T|z_1(s,t)-z_2(s,t)|ds\)^2\]dt\\
\ns\ds\le K\dbE\Big\{\int_0^Te^{\b
t}\[|\psi_1(t)-\psi_2(t)|^2+\(\int_t^T|\d
g(t,s)|ds\)^2\]dt\\
\ns\ds\qq+{1\over\b}\|(y_1(\cd),z_1(\cd\,,\cd))-(y_2(\cd),z_2(\cd\,,\cd))\|^2_{\cM^2_\b[0,T]}
\Big\}.\ea\ee
Hence, by letting $\psi_i(\cd)=\psi(\cd)$ and $g_i(\cd)=g(\cd)$, we
have
$$\|\Th(y_1(\cd),z_1(\cd\,,\cd))-\Th(y_2(\cd),z_2(\cd\,,\cd))\|_{\cM^2_\b[0,T]}^2
\le{K\over\b}\|(y_1(\cd),z_1(\cd\,,\cd))-(y_2(\cd),z_2(\cd\,,\cd))\|^2_{\cM^2_\b[0,T]}.$$
This implies that $\Th$ can be naturally extended to $\cM^2[0,T]$
since $\wt\cM^2[0,T]$ is dense in $\cM^2[0,T]$. Moreover, since the
constant $K>0$ appears in the right hand side of the above is
independent of $\b$, by choosing $\b>0$ large enough, we obtain that
the extension of $\Th$, still denoted by itself, is a contraction.
Hence, $\Th$ admits a unique fixed point
$(Y(\cd),Z(\cd\,,\cd))\in\cM^2[0,T]$, which is the unique adapted
M-solution of (\ref{2.3}).

\ms

Note that (\ref{2.18}) implies that for general $\psi(\cd)\in
L^2_{\cF_T}(0,T;\dbR^m)$ and general generator $g(\cd)$ satisfying
(H1), by an approximating argument, one can obtain that the
corresponding BSVIE admits a unique adapted M-solution
$(Y(\cd),Z(\cd\,,\cd))$. Also, the stability estimate
(\ref{stability}) follows.

\ms

On the other hand, similar to (\ref{2.18}), we have
$$\ba{ll}
\ns\ds\|(Y(\cd),Z(\cd\,,\cd))\|^2_{\cM^2_\b[0,T]}\equiv
\int_0^Te^{\b t}\dbE\[|Y(t)|^2+\int_t^T|Z(t,s)|^2ds\]dt\\
\ns\ds\le K\int_0^Te^{\b
t}\dbE\[|\psi(t)|^2+\(\int_t^T|g_0(t,s)|ds\)^2+\int_t^T|Y(s)|^2ds+\(\int_t^T|Z(s,t)|ds\)^2\]dt\\
\ns\ds\le K\dbE\Big\{\int_0^Te^{\b
t}\[|\psi(t)|^2+\(\int_t^T|g_0(t,s)|ds\)^2\]dt+{1\over\b}\|(Y(\cd),Z(\cd\,,\cd))\|^2_{\cM^2_\b[0,T]}
\Big\},\ea$$
which leads to the estimate (\ref{2.7}).

\ms

(ii) Let $(Y(\cd),Z(\cd\,,\cd))$ be the adapted solution of BSVIE
(\ref{2.4}). For any fixed $t\in[0,T]$, we let
$(\eta(t,\cd),\z(t,\cd))$ be the adapted solution to the following
BSDE:
\bel{BSDE2}\eta(t,r)=\psi(t)+\int_r^Tg(t,s,Y(s),\z(t,s))ds-\int_r^T\z(t,s)dW(s),\qq
r\in[0,T].\ee
Then we know that
$$Y(t)=\eta(t,t),\qq Z(t,s)=\z(t,s),\qq0\le t\le s\le T.$$
By (\ref{BSDE2}), we have
\bel{2.20}\ba{ll}
\ns\ds\dbE\[\sup_{r\in[t,T]}|\eta(t,r)|^2+\int_t^T|\z(t,s)|^2ds\]\le
K\dbE\[|\psi(t)|^2+\(\int_t^T|g(t,s,Y(s),0)|ds\)^2\]\\
\ns\ds\le
K\dbE\[|\psi(t)|^2+\(\int_0^T|g_0(t,s)|ds\)^2+\int_t^T|Y(s)|^2ds\].\ea\ee
Thus,
$$\ba{ll}
\ns\ds\dbE\[|Y(t)|^2+\int_t^T|Z(t,s)|^2ds\]\le
K\dbE\[|\psi(t)|^2+\(\int_0^T|g_0(t,s)|ds\)^2+\int_t^T|Y(s)|^2ds\].\ea$$
By Gronwall's inequality, we obtain estimate (\ref{2.9}). This also
leads to
\bel{2.22}\dbE\[\sup_{r\in[t,T]}|\eta(t,r)|^2+\int_t^T|\z(t,s)|^2ds\]
\le K\dbE\[|\psi(t)|^2+\(\int_0^T|g_0(t,s)|ds\)^2\].\ee
Similar to (\ref{2.16}), in the current case, we have
$$\ba{ll}
\ns\ds\dbE\[|Y_1(t)-Y_2(t)|^2+\int_t^T|Z_1(t,s)-Z_2(t,s)|^2ds\]\\
\ns\ds\le K\dbE\[|\psi_1(t)-\psi_2(t)|^2+\(\int_t^T|\d
g(t,s)|ds\)^2+\int_t^T|Y_1(s)-Y_2(s)|^2ds\].\ea$$
Then applying Gronwall's inequality, we obtain stability estimate
(\ref{2.10}). To prove the continuity of $t\mapsto Y(t)$, we let
$t,t'\in[0,T]$ and consider the following:
$$\ba{ll}
\ns\ds\eta(t,r)-\eta(t',r)=\psi(t)-\psi(t')+\int_r^T\(g(t,s,Y(s),\z(t,s))-g(t',s,Y(s),
\z(t',s))\)ds\\
\ns\ds\qq\qq\qq\qq-\int_r^T\(\z(t,s)-\z(t',s)\)dW(s),\qq
r\in[0,T].\ea$$
Then the stability of adapted solutions to BSDEs implies that
$$\ba{ll}
\ns\ds\dbE\[\sup_{r\in[0,T]}|\eta(t,r)-\eta(t',r)|^2+\int_0^T
|\z(t,s)-\z(t',s)|^2ds\]\\
\ns\ds\le
K\dbE\[|\psi(t)-\psi(t')|^2+\(\int_0^T|g(t,s,Y(s),0)-g(t',s,Y(s),0)|ds\)^2\]\\
\ns\ds\le
K\dbE\[|\psi(t)-\psi(t')|^2+\rho(|t-t'|)^2\(\int_0^T|Y(s)|^2ds\)\]\\
\ns\ds\le
K\dbE\Big\{|\psi(t)-\psi(t')|^2+\rho(|t-t'|)^2\[\int_0^T|\psi(t)|^2dt+\int_0^T\(\int_0^T|g(t,s,0,0)|
ds\)^2dt\]\Big\}.\ea$$
Therefore, one has
$$\left\{\2n\ba{ll}
\ns\ds\lim_{|t-t'|\to0}\dbE|\eta(t,r)-\eta(t',r)|^2=0,\qq\hb{uniformly
in $r\in[0,T]$},\\
\ns\ds\lim_{r'\to r}\dbE|\eta(t,r)-\eta(t,r')|^2=0,\qq\forall
(t,r)\in[0,T]^2.\ea\right.$$
Hence, $(t,r)\mapsto\eta(t,r)$ is continuous, i.e.,
$$\lim_{(t',r')\to(t,r)}\dbE|\eta(t',r')-\eta(t,r)|^2=0,\qq\forall(t,r)\in[0,T]^2.$$
Consequently, $t\mapsto\eta(t,t)=Y(t)$ is continuous from $[0,T]$ to
$L^2_{\cF_T}(\O;\dbR^m)$, which leads to $Y(\cd)\in C_{\dbF}(0,T;$
$L^2(\Omega;\dbR^m))$. \endpf

\section{Optimal Control Problems and Maximum Principles}

Now, we consider the following controlled FBSVIE:
\bel{FBSVIE3.1}\left\{\2n\1n\ba{ll}
\ns\ds
X(t)=\f(t)+\int_0^tb(t,s,X(s),u(s))ds+\int_0^t\si(t,s,X(s),u(s))dW(s), \\
\ns\ds
Y(t)=\psi(t,X(t),X(T))+\int_t^Tg(t,s,X(t),X(s),Y(s),Z(t,s),Z(s,t),u(s))ds\\
\ns\ds\qq\qq\q-\int_t^TZ(t,s)dW(s),\qq t\in[0,T],\ea\right.\ee
where admissible control $u(\cd)$ belongs to $\cU[0,T]$ defined by
$$\cU[0,T]=\Big\{u(\cdot )\in L_{\mathcal{F}}^2(0,T;\dbR^\ell)\bigm|u(t)\in
U,~\as\2n,~\ae t\in[0,T]\Big\},$$
with $U$ being a nonempty convex subset of $\dbR^\ell$. For
convenience, we let $0\in U$. Also, we will consider the following
FBSVIE which is a special case of (\ref{FBSVIE3.1}):
\bel{FBSVIE3.2}\left\{\2n\1n\ba{ll}
\ns\ds
X(t)=\f(t)+\int_0^tb(t,s,X(s),u(s))ds+\int_0^t\si(t,s,X(s),u(s))dW(s), \\
\ns\ds
Y(t)=\psi(t,X(t),X(T))+\int_t^Tg(t,s,X(t),X(s),Y(s),Z(t,s),u(s))ds\\
\ns\ds\qq\qq\q-\int_t^TZ(t,s)dW(s),\qq t\in[0,T],\ea\right.\ee
A triple
$$(X(\cd),Y(\cd),Z(\cd\,,\cd))\in
C_{\dbF}(0,T;L^2(\Omega;\dbR^n))\times L_\dbF^2(0,T;\dbR^m)\times
L^2(0,T;L_\dbF^2(0,T;\dbR^m))$$
is called an adapted M-solution of (\ref{FBSVIE3.1}) if $X(\cd)$
satisfies the forward stochastic Volterra integral equation (FSVIE,
for short) in (\ref{FBSVIE3.1}) and $(Y(\cd),Z(\cd,\cd))$ is the
adapted M-solution of the BSVIE in (\ref{FBSVIE3.1}). Also, a
triple
$$(X(\cd),Y(\cd),Z(\cd))\in C_{\dbF}(0,T;L^2(\Omega;\dbR^n))\times C_\dbF([0,T];L^2(\O;\dbR^m))\times L^2_\dbF(\D;
\dbR^m)$$
is called an adapted solution of (\ref{FBSVIE3.2}) if $X(\cd)$
satisfies the FSVIE in (\ref{FBSVIE3.2}) and $(Y(\cd),Z(\cd,\cd))$
is an adapted solution of BSVIE in (\ref{FBSVIE3.2}).

\ms

The following collects some basic assumptions on the coefficients of
FBSVIE (\ref{FBSVIE3.1}).

\ms

{\bf(H3)} Let $\varphi(\cd)\in C_{\dbF}(0,T;L^2(\Omega;\dbR^n))$,
$\psi(\cd,0,0)\in L^2_{\cF_T}(0,T;\dbR^m)$, and let
$$\left\{\2n\ba{ll}
\ns\ds b,\si:[0,T]^2\times\dbR^n\times U\times\O\to\dbR^n,\\
\ns\ds
g:[0,T]^2\times\dbR^n\times\dbR^n\times\dbR^m\times\dbR^m\times\dbR^m\times
U\times\O\to\dbR^m\ea\right.$$
be measurable such that
$$s\mapsto\big(b(t,s,x,u),\si(t,s,x,u),g(t,s,x',x,y,z,z',u)\big)$$
is $\dbF$-progressively measurable on $[0,T]$,
$$(x',x,y,z,z',u)\mapsto\big(b(t,s,x,u),\si(t,s,x,u),
g(t,s,x',x,y,z,z',u),\psi(t,x',x)\big)$$
is continuously differentiable with uniformly bounded derivatives,
and with the notation
$$b_0(t,s)\equiv b(t,s,0,0),\q\si_0(t,s)\equiv\si(t,x,0,0),\q g_0(t,s)\equiv g(t,s,0,0,0,0,0,0),$$
one has
\bel{}\sup_{t\in[0,T]}\dbE\[\(\int_0^T|b_0(t,s)|ds\)^2+\int_0^T|\si_0(t,s)|ds\]
+\dbE\int_0^T\(\int_0^T|g_0(t,s)|ds\)^2dt <\infty.\ee
Further, there exists a modulus of continuity
$\rho:[0,\infty)\to[0,\infty)$ (i.e., $\rho(\cd)$ is continuous and
strictly increasing with $\rho(0)=0$) such that
$$\ba{ll}
\ns\ds|b(t,s,x,u)-b(t',s,x,u)|+|\si(t,s,x,u)-\si(t',s,x,u)|\le\rho(|t-t'|)(1+|x|+|u|),\\
\ns\ds\qq\qq\qq\qq\qq\qq\qq\qq\qq
t,t',s\in[0,T],~(x,u)\in\dbR^n\times U.\ea$$

\ms

For FBSVIE (\ref{FBSVIE3.2}), we introduce the following stronger
hypothesis.

\ms

{\bf(H4)} Let $\varphi(\cd)\in C_{\dbF}(0,T;L^2(\Omega;\dbR^n))$,
$\psi(\cd,0,0)\in C_{\dbF}(0,T;L^2(\Omega;\dbR^m))$, and let
$$\left\{\2n\ba{ll}
\ns\ds b,\si:[0,T]^2\times\dbR^n\times U\times\O\to\dbR^n,\\
\ns\ds
g:[0,T]^2\times\dbR^n\times\dbR^n\times\dbR^m\times\dbR^m\times
U\times\O\to\dbR^m,\ea\right.$$
be measurable such that
$$s\mapsto\big(b(t,s,x,u),\si(t,s,x,u)),g(t,s,x',x,y,z,u)\big)$$
is $\dbF$-progressively measurable on $[0,T]$,
$$(x',x,y,z,u)\mapsto(b(t,s,x,u),\si(t,s,x,u),g(t,s,x',x,y,z,u),\psi(t,x',x))$$
is continuously differentiable with uniformly bounded derivatives,
and with notation
$$b_0(t,s)\equiv b(t,s,0,0),\q\si_0(t,s)\equiv\si(t,x,0,0),\q g_0(t,s)\equiv g(t,s,0,0,0,0,0),$$
one has
\bel{}\sup_{t\in[0,T]}\dbE\[\(\int_0^T|b_0(t,s)|ds\)^2+\int_0^T|\si_0(t,s)|ds+\int_0^T|g_0(t,s)|ds\)^2\]
<\infty.\ee
Further, there exists a modulus of continuity
$\rho:[0,\infty)\to[0,\infty)$ such that
$$\ba{ll}
\ns\ds|b(t,s,x,u)-b(t',s,x,u)|+|\si(t,s,x,u)-\si(t',s,x,u)|+|\psi(t,x',x)-\psi(t',x',x)|\\
\ns\ds\q+|g(t,s,x',x,y,z,u)-g(t',s,x',x,y,z,u)|\le\rho(|t-t'|)(1+|x|+|x'|+|y|+|z|+|u|),\\
\ns\ds\qq\qq\qq\qq\qq\qq\qq
t,t',s\in[0,T],~x,x'\in\dbR^n,~y,z\in\dbR^m,~u\in U.\ea$$

\ms

The following result follows from some standard theory of stochastic
Volterra integral equations and those presented in the previous
section.

\ms

\bf Proposition 3.1. \sl {\rm(i)} Let {\rm(H3)} hold. Then for any
$u(\cd)\in\cU[0,T]$, (\ref{FBSVIE3.1}) admits a unique adapted
M-solution $(X(\cd),Y(\cd),Z(\cd\,,\cd))$. Moreover, the following
holds:
\bel{}\ba{ll}
\ns\ds\sup_{t\in[0,T]}\dbE|X(t)|^2+\dbE\int_0^T|Y(t)|^2dt+\dbE\int_0^T\int_0^T|Z(t,s)|^2dsdt\\
\ns\ds\le
K\Big\{\sup_{t\in[0,T]}\dbE|\f(t)|^2+\sup_{t\in[0,T]}\dbE\(\int_0^T|b_0(t,s)|ds\)^2
+\sup_{t\in[0,T]}\dbE\int_0^T|\si_0(t,s)|^2ds\\
\ns\ds\qq+\dbE\int_0^T|\psi(t,0,0)|^2dt+\dbE\int_0^T\int_0^T|g_0(t,s)|ds\)^2\Big\}.\ea\ee
In addition, if
$(\f_i(\cd),b_i(\cd),\si_i(\cd),\psi_i(\cd),g_i(\cd))$ also satisfy
$(H3)$, $u_i(\cd)\in\cU[0,T]$, $i=1,2$, and $(X_i(\cd),Y_i(\cd)$,
$Z_i(\cd\,,\cd))$ is the adapted M-solution of the corresponding
FBSVIEs, then
\bel{}\ba{ll}
\ns\ds\sup_{t\in[0,T]}\dbE|X_1(t)-X_2(t)|^2dt+\dbE\2n\int_0^T\2n|Y_1(t)-Y_2(t)|^2dt+\dbE\2n
\int_0^T\2n\int_0^T
\2n|Z_1(t,s)-Z_2(t,s)|^2dsdt\\
\ns\ds\le
K\Big\{\sup_{t\in[0,T]}\dbE|\f_1(t)\1n-\1n\f_2(t)|^2+\2n\sup_{t\in[0,T]}\dbE\(\int_0^T\2n|
b_1(t,s,X_1(s),u_1(s))\1n-\1n b_2(t,s,X_1(s),
u_2(s))|ds\)^2\\
\ns\ds\qq+\sup_{t\in[0,T]}\dbE\int_0^T|\si_1(t,s,X_1(s),u_1(s))-\si_2(t,s,X_1(s),u_2(s))|^2ds\\
\ns\ds\qq+\dbE\int_0^T\(\int_0^T|g_1(t,s,X_1(t),X_1(s),Y_1(s),Z_1(t,s),Z_1(s,t),u_1(s))\\
\ns\ds\qq\qq\qq\qq-g_2(t,s,X_1(t),X_1(s),Y_1(s),Z_1(t,s),Z_1(s,t),u_2(s))|ds\)^2dt\Big\}.\ea\ee

\ms

{\rm(ii)} Let {\rm(H4)} hold. Then for any $u(\cd)\in\cU[0,T]$,
(\ref{FBSVIE3.2}) admits a unique adapted solution
$(X(\cd),Y(\cd),Z(\cd\,,\cd))$. Moreover, the following holds:
\bel{}\ba{ll}
\ns\ds\sup_{t\in[0,T]}\dbE|X(t)|^2+\sup_{t\in[0,T]}\dbE|Y(t)|^2+\sup_{t\in[0,T]}\dbE\int_t^T|Z(t,s)|^2ds\\
\ns\ds\le
K\Big\{\sup_{t\in[0,T]}\dbE|\f(t)|^2+\sup_{t\in[0,T]}\dbE\(\int_0^T|b_0(t,s)|ds\)^2
+\sup_{t\in[0,T]}\dbE\int_0^T|\si_0(t,s)|^2ds\\
\ns\ds\qq+\sup_{t\in[0,T]}\dbE|\psi(t,0,0)|^2+\sup_{t\in[0,T]}\dbE\(\int_0^T|g_0(t,s)|ds\)^2\Big\}.\ea\ee
In addition, if
$(\f_i(\cd),b_i(\cd),\si_i(\cd),\psi_i(\cd),g_i(\cd))$ satisfy
$(H4)$, $u_i(\cd)\in\cU[0,T]$, $i=1,2$, and $(X_i(\cd),Y_i(\cd)$,
$Z_i(\cd\,,\cd))$ is the adapted M-solution of the corresponding
FBSVIEs, then
\bel{}\ba{ll}
\ns\ds\sup_{t\in[0,T]}\dbE|X_1(t)\1n-\1n
X_2(t)|^2dt\1n+\2n\sup_{t\in[0,T]}\dbE|Y_1(t)\1n-\1n
Y_2(t)|^2dt\1n+\2n\sup_{t\in[0,T]}\dbE\2n\int_t^T
\2n|Z_1(t,s)\1n-\1n Z_2(t,s)|^2dsdt\\
\ns\ds\le\1n
K\Big\{\sup_{t\in[0,T]}\dbE|\f_1(t)\1n-\1n\f_2(t)|^2\2n+\2n\sup_{t\in[0,T]}\dbE\(\int_0^T\2n|
b_1(t,s,X_1(s),u_1(s))\1n-\1n b_2(t,s,X_1(s),u_2(s))|ds\)^2\\
\ns\ds\qq+\sup_{t\in[0,T]}\dbE\int_0^T|\si_1(t,s,X_1(s),u_1(s))-\si_2(t,s,X_1(s),u_2(s))|^2ds\\
\ns\ds\qq+\sup_{t\in[0,T]}\dbE\(\int_0^T|g_1(t,s,X_1(t),X_1(s),Y_1(s),Z_1(t,s),u_1(s))\\
\ns\ds\qq\qq\qq\qq-g_2(t,s,X_1(t),X_1(s),Y_1(s),Z_1(t,s),u_2(s))|ds\)^2\Big\}.\ea\ee

\ms

\rm

We see that whether the generator $g(t,s,y,z,z')$ depends on $z'$
will have different regularity of $Y(\cd)$ in general. Therefore, we
will introduce two different optimal control problems.

\ms

First, we consider state equation (\ref{FBSVIE3.2}). Since for such
a case, $Y(0)$ is well defined, we may introduce the cost functional
as follows:
\bel{cost3.1}J_1(u(\cd))=
\dbE\[h(X(T),Y(0))+\int_0^T\int_t^Tf(t,s,X(s),Y(s),Z(t,s),u(s))dsdt\].\ee
For the involved functions $h$ and $f$ in (\ref{cost3.1}), we impose
the following hypothesis.

\ms

{\bf (H5)} Let
$$\ba{ll}
\ns\ds h:\dbR^n\times\dbR^m\times\O\to\dbR,\qq
f:\Delta\times\dbR^n\times\dbR^m\times\dbR^m \times U \times\O\to\dbR\ea$$
be measurable such that
$$(x,y)\mapsto h(x,y),\q(x,y,z,u)\mapsto f(t,s,x,y,z,u)$$
are continuously differentiable with the derivatives of $h$ and $f$
being bounded by $K(1+|x|+|y|)$ and $K(1+|x|+|y|+|z|+|u|)$,
respectively.

\ms

Now, we state our first optimal control problem.

\ms

\bf Problem (C1). \rm With the state equation (\ref{FBSVIE3.2}), find
$\bar u(\cd)$ such that
\bel{3.4}J_1(\bar u(\cd))=\inf_{u(\cd)\in\cU[0,T]}J_1(u(\cd)).\ee

\ms

Any $\bar u(\cd)\in\cU[0,T]$ satisfying (\ref{3.4}) is called an
{\it optimal control} of Problem (C1). The corresponding state
process, denoted by $(\bar X(\cd),\bar Y(\cd),\bar Z(\cd\,,\cd))$,
is called an {\it optimal state process}, and $(\bar X(\cd),\bar
Y(\cd),\bar Z(\cd\,,\cd),\bar u(\cd))$ is called an {\it optimal
4-tuple} of Problem (C1). To make the statement of the maximum
principle for Problem (C1) simpler, let us introduce the following
notations which will also be used in Section 5. For any given
optimal 4-tuple $(\bar X(\cd),\bar Y(\cd),\bar Z(\cd\,,\cd),\bar
u(\cd))$ of Problem (C1), we denote
$$b_x(t,s)=b_x(t,s,\bar X(s),\bar u(s)),\q b_u(t,s)=b(t,s,\bar
X(s),\bar u(s)).$$
The notations $\si_x(t,s)$, $\si_u(t,s)$, $\psi_{x'}(t)$,
$\psi_x(t)$, $g_{x'}(t,s)$, $g_x(t,s)$, $g_y(t,s)$, $g_z(t,s)$,
$g_{z'}(t,s)$, $g_u(t,s)$, $h_x$, $h_y$, $f_x(t,s)$, $f_y(t,s)$, and
$f_u(t,s)$ are similar. Also, for any scalar valued function, say
$x\mapsto f(t,s,x,y,z,u)$, $f_x(t,s,x,y,z,u)$ is regarded as row
vector, i.e., $\dbR^{1\times n}$-valued. Such a convention will be
consistent with vector valued functions, say, $x\mapsto\psi(t,x,x')$
for which $\psi_x(t,x,x')$ takes values in $\dbR^{m\times n}$. We
now state the following maximum principle.

\ms

\bf Theorem 3.2. \sl Let {\rm(H4)} and {\rm(H5)} hold. Let $(\bar
X(\cd),\bar Y(\cd),\bar Z(\cd\,,\cd),\bar u(\cd))$ be an optimal
4-tuple of Problem {\rm(C1)}. Then
\bel{maximum1}\ba{ll}
\ns\ds\lan\int_s^T\big[b_u(t,s)^Tp(t)+\si_u(t,s)^T
q(t,s)\big]dt+b_u(T,s)^T\m(T)+\si_u(T,s)^T\n(s)\\
\ns\ds\qq+g_u(0,s)^T\l(s)+\int_0^sf_u(t,s)^Tdt
+\int_0^sg_u(t,s)^T\xi(t)dt,v-\bar u(s)\ran\ge0,\\
\ns\ds\qq\qq\qq\qq\qq\qq\qq\forall v\in U,~\ae t\in[0,T],~\as,\ea\ee
where
\bel{adjoint1}\left\{\2n\ba{ll}
\ns\ds\l(t)=\dbE(h_y^T)+\int_0^tg_z(0,s)^T\l(s)dW(s),\\
\ns\ds\xi(t)=g_y(0,t)^T\l(t) +\int_0^tf_y(s,t)^Tds
+\int_t^Tf_z(t,s)^TdW(s)\\
\ns\ds\qq\qq+\int_0^tg_y(s,t)^T\dbE_t[\xi(s)]ds+\int_t^Tg_z(t,s)^T\dbE_s[\xi(t)]dW(s),\\
\ns\ds\m(t)=h_x^T+\psi_x(0)^T\l(T)+\int_0^T\psi_x(s)^T\xi(s)ds
-\int_t^T\n(s)dW(s),
\\
\ns\ds p(t)=b_x(T,t)^T\m(T)+\si_x(T,t)^T\n(t)+g_x(0,t)^T\l(t)
+\int_0^tf_x(s,t)^Tds\\
\ns\ds\qq\qq+\(\psi_{x'}(t)^T
+\int_t^Tg_{x'}(t,s)^Tds\)\xi(t)+\int_0^tg_x(s,t)^T\xi(s)ds\\
\ns\ds\qq\qq+\int_t^T\(b_u(s,t)^Tp(s)+\si_u(s,t)^Tq(s,t)\)ds
-\int_t^Tq(t,s)dW(s).\ea\right.\ee

\ms

\rm

Note that in (\ref{adjoint1}), $\l(\cd)$ solves an FSDE, $\xi(\cd)$
solves a special type of stochastic Fredholm integral equation with
mean-field. We will show in the next section that such an equation
admits a unique solution $\xi(\cd)$ which is not required to be
$\dbF$-adapted. The equation for $(\m(\cd),\n(\cd))$ is a BSDE, and
that for $(p(\cd),q(\cd\,,\cd))$ is a BSVIE. We see that the system
(\ref{adjoint1}) is a decoupled system.

\ms

Next let us consider the state equation (\ref{FBSVIE3.1}). For this
case, we introduce the following cost functional:
\bel{cost3.2}J_2(u(\cd))=
\dbE\[h\left(X(T),\dbE\int_0^TY(s)ds\right)
+\int_0^T\int_0^Tf(t,s,X(s),Y(s),Z(t,s),u(s))dsdt\].\ee
For the involved functions $h$ and $f$ in (\ref{cost3.2}), we impose
the following hypothesis.

\ms

{\bf (H6)} Let
$$\ba{ll}
\ns\ds h:\dbR^n\times\dbR^m\times\O\to\dbR,\qq
f:[0,T]^2\times\dbR^n\times\dbR^m\times\dbR^m \times U \times\O\to\dbR\ea$$
be measurable such that
$$(x,y)\mapsto h(x,y),\q(x,y,z,u)\mapsto f(t,s,x,y,z,u)$$
are continuously differentiable with the derivatives of $h$ and $f$
being bounded by $K(1+|x|+|y|)$ and $K(1+|x|+|y|+|z|+|u|)$,
respectively.

\ms

We may pose the following problem.

\ms

\bf Problem (C2). \rm With the state equation (\ref{FBSVIE3.1}), find
$\bar u(\cd)$ such that
\bel{3.10}J_2(\bar u(\cd))=\inf_{u(\cd)\in\cU[0,T]}J_2(u(\cd)).\ee

\ms

Any $\bar u(\cd)\in\cU[0,T]$ satisfying (\ref{3.10}) is called an
optimal control of Problem (C2). The corresponding state process,
denoted by $(\bar X(\cd),\bar Y(\cd),\bar Z(\cd\,,\cd))$, is called
an optimal state process, and $(\bar X(\cd),\bar Y(\cd),\bar
Z(\cd\,,\cd),\bar u(\cd))$ is called an optimal 4-tuple of Problem
(C2). We have the following maximum principle.

\ms

\bf Theorem 3.3. \sl Let {\rm(H3)} and {\rm(H6)} hold. Let $(\bar
X(\cd),\bar Y(\cd),\bar Z(\cd\,,\cd),\bar u(\cd))$ be the
corresponding M-solution of FBSVIE (\ref{FBSVIE3.1}). Then
$$\ba{ll}
\ns\ds\lan b_u(T,t)^T\m(T)+\si_u(T,t)^T\n(t)+
\int_0^Tf_u(s,t)^Tds+\int_0^tg_u(s,t)^T\xi(s)ds\\
\ns\ds\qq\qq+\int_t^Tb_u(s,t)^Tp(s)ds+\int_t^T\si_u(s,t)^Tq(s,t)ds,
v-\bar u(t)\ran\ge0,\\
\ns\ds\qq\qq\qq\qq\qq\qq\qq\qq\qq\forall v\in U,~\as\2n,~\ae
t\in[0,T],\ea$$
where $(\xi(\cd),\m(\cd),\n(\cd),p(\cd),g(\cd\,,\cd))$ solves the
following adjoint equation:
\bel{adjoint2}\left\{\ba{ll}
\ns\ds\xi(t)=h_y^T\1n+\2n\int_0^Tf_y(s,t)^Tds+\2n\int_0^Tf_z(t,s)^TdW(s)\\
\ns\ds\qq\qq+\int_0^tg_y(s,t)^T\dbE_t[\xi(s)]ds
+\int_0^t\dbE_s[g_{z'}(s,t)^T\xi(s)]dW(s)\\
\ns\ds\qq\qq+\int_t^Tg_z(t,s)^T\dbE_s[\xi(t)]dW(s),\\
\ns\ds\m(t)=h_x^T+\int_0^T\psi_x(s)^T\xi(s)ds-\int_t^T\n(s)dW(s),\\
\ns\ds
p(t)=b_x(T,s)^T\m(T)+\si_x(T,s)^T\n(s)+\int_0^Tf_x(s,s)^Tds\\
\ns\ds\qq\qq+\(\psi_{x'}(s)^T
+\int_t^Tg_{x'}(t,s)^Tds\)\xi(t)+\int_0^sg_x(s,t)^T\xi(s)
ds\\
\ns\ds\qq\qq+\int_t^T\[b_x(s,t)^Tp(s)+\si_x(s,t)^Tq(s,t)\]ds-\int_t^Tq(t,s)
dW(s).\ea\right.\ee

\ms

\rm

We see that different from Theorem 3.2, in the above the adjoint
equation only consists of three equations, the equation for
$\l(\cd)$ is not necessary here. Actually, we will see that the
equation for $\l(\cd)$ is used to take care of the term involving
$Y(0)$ in the cost functional. Again, (\ref{adjoint2}) is also
decoupled.

\section{Duality Principles}

\rm

The aim of this section is to establish a duality principle between
the following linear BSVIE:
\bel{4.1}\ba{ll}
\ns\ds
Y(t)=\psi(t)+\int_t^T\big[A(t,s)Y(s)+B(t,s)Z(t,s)+C(t,s)Z(s,t)\big]ds\\
\ns\ds\qq\qq\qq\qq\qq-\int_t^TZ(t,s)dW(s),\qq t\in[0,T],\ea\ee
and its adjoint equation, where $(Y(\cd),Z(\cd,\cd))$ is the unique
M-solution associated with $\psi(\cd)\in L^2_{\cF_T}(0,T;\dbR^m)$.
We introduce the following hypothesis for the coefficients of the
above equation.

\ms

{\bf(H7)} The maps $A,B,C:[0,T]^2\times\O\to\dbR^{m\times m}$ are
measurable and bounded such $s\1n\mapsto\1n(A(t,s),\1n B(t,s),\1n
C(t,s))$ is $\dbF$-adapted on $[t,T]$ for every $t\in[0,T]$.

\ms

By Theorem 2.3, we know that under (H7), for any $\psi(\cd)\in
L^2_{\cF_T}(0,T;\dbR^m)$, linear BSVIE (\ref{4.1}) admits a unique
adapted M-solution $(Y(\cd),Z(\cd\,,\cd))\in\cM^2[0,T]$. Note that
in \cite{Yong 2008} (see Theorem 5.1 there), a duality principle was
established for the case that $Z(t,s)$ does not appear (or
$B(\cd\,,\cd)=0$). The significance here in the current paper is
that we have discovered the adjoint equation of (\ref{4.1}) with all
the interested terms appearing and we have the well-posedness of
such an equation. We now introduce the adjoint equation for
(\ref{4.1}). For any $(\a(\cd),\b(\cd\,,\cd))\in
L^2_\dbF(0,T;\dbR^m)\times L^2(0,T;L^2_{\dbF}(0,T;\dbR^m))$,
consider the following stochastic integral equation:
\bel{4.2}\ba{ll}
\ns\ds\xi(t)=\a(t)+\int_0^tA(s,t)^T\dbE_t[\xi(s)]ds+\int_0^t\dbE_s[C(s,t)^T\xi(s)]dW(s)\\
\ns\ds\qq\qq+\int_0^T\b(t,s)dW(s)+\int_t^TB(t,s)^T\dbE_s[\xi(t)]dW(s),
\qq\as\2n,~\ae t\in[0,T],\ea\ee
where $\dbE_r[\z]\deq\dbE[\z\,|\,\cF_r]$ for any integrable random
variable $\z$ and $r\in[0,T]$. We call (\ref{4.2}) the {\it adjoint
equation} of linear BSVIE (\ref{4.1}). It is seen that (\ref{4.2})
is a mean-field stochastic Fredholm-Volterra type integral equation
with some special structure, whose unknown is an $\cF_T$-measurable
process $\xi(\cd)$. Unlike usual BSDEs or BSVIEs, in the above, we
do not require $\xi(\cd)$ to be $\dbF$-adapted. We now state the
duality principle.

\ms

\bf Theorem 4.1. \sl Let {\rm(H7)} hold. Then for any
$(\a(\cd),\b(\cd\,,\cd))\in L^2_\dbF(0,T;\dbR^n)\times
L^2(0,T;L^2_{\dbF}(0,T;\dbR^n))$, equation $(\ref{4.2})$ admits a
unique solution $\xi(\cd)\in L^2_{\cF_T}(0,T;\dbR^m)$. Further, if
$(Y(\cd),Z(\cd,\cd))$ is the unique M-solution for $(\ref{4.1})$
with $\psi(\cd)\in L^2_{\cF_T}(0,T;\dbR^n)$, then
\bel{duality}\ba{ll}
\ns\ds\dbE\int_0^T\lan\psi(t),\xi(t)\ran dt=\dbE\Big\{\int_0^T\lan
Y(t),\a(t)\ran dt+\int_0^T\int_0^T\lan Z(t,s),\b(t,s)\ran
dsdt\Big\}.\ea\ee

\ms

\it Proof. \rm We first prove the well-posedness of the adjoint
equation (\ref{4.2}). Let $(\a(\cd),\b(\cd\,,\cd))\in
L^2_\dbF(0,T;\dbR^n)\times L^2(0,T;L^2_{\dbF}(0,T;\dbR^n))$ be
given. Suppose $\xi(\cd)$ is a solution to (\ref{4.2}). For almost
all $t\in[0,T]$ and any $r\in[t,T]$, applying $\dbE_r$ on the both
sides of (\ref{4.2}), one gets
\bel{4.4}\ba{ll}
\ns\ds\dbE_r[\xi(t)]=\a(t)+\int_0^r\b(t,s)dW(s)+\int_t^rB(t,s)^T\dbE_s[\xi(t)]dW(s)\\
\ns\ds\qq\qq\qq+\int_0^tA(s,t)^T\dbE_t[\xi(s)]ds+\int_0^t\dbE_s[C(s,t)^T\xi(s)]dW(s)\\
\ns\ds\qq\qq=\a(t)+\int_0^t\b(t,s)dW(s)+\int_t^r\(B(t,s)^T\dbE_s[\xi(t)]+\b(t,s)\)dW(s)\\
\ns\ds\qq\qq\qq+\int_0^tA(s,t)^T\dbE_t[\xi(s)]ds\1n+\1n\int_0^t\dbE_s[C(s,t)^T\xi(s)]dW(s),
\q \forall r\in[t,T],~t\in[0,T].\ea\ee
On the other hand, if $\xi(\cd)$ satisfies (\ref{4.4}) for any
$r\in[t,T]$, then taking $r=T$, we recover (\ref{4.2}). Thus, to
solve (\ref{4.2}), it suffices to solve (\ref{4.4}).

\ms

For any $\bar\xi(\cd)\in L^2_{\cF_T}(0,T;\dbR^m)$, and almost every
$t\in[0,T]$, consider the following family of stochastic
differential equation with parameter $t$:
\bel{4.5}\ba{ll}
\ns\ds\l(t,r)=\a(t)+\int_0^t\b(t,s)dW(s)+\int_t^r\(B(t,s)^T\l(t,s)+\b(t,s)\)dW(s)\\
\ns\ds\qq\qq\qq+\int_0^tA(s,t)^T\dbE_t[\bar\xi(s)]ds+\int_0^t\dbE_s[C(s,t)^T
\bar\xi(s)]dW(s), \qq r\in[t,T].\ea\ee
It admits a unique solution $\l(t,\cd)$ which is $\dbF$-adapted on
$[t,T]$, and
$$\dbE|\l(t,r)|^2\le
K\dbE\Big\{|\a(t)|^2+\int_0^r|\b(t,s)|^2ds+\int_t^r|\l(t,s)|^2ds
+\int_0^t|\bar\xi(s)|^2ds\Big\},\qq r\in[t,T].$$
Hence, it follows from Gronwall's inequality that
\bel{4.6}\dbE|\l(t,r)|^2\le
K\dbE\Big\{|\a(t)|^2+\int_0^r|\b(t,s)|^2ds+\int_0^t|\bar\xi(s)|^2ds\Big\},\qq
r\in[t,T].\ee
We define
$$\xi(t)=\l(t,T),\qq\ae t\in[0,T].$$
Then for any $r\in[t,T]$, from (\ref{4.5}), one has
\bel{}\dbE_r[\xi(t)]=\dbE_r[\l(t,T)]=\l(t,r),\ee
and it follows from (\ref{4.6}) that
$$\dbE|\xi(t)|^2\le K\dbE\Big\{|\a(t)|^2+\int_0^T|\b(t,s)|^2ds
+\int_0^t|\bar\xi(s)|^2ds\Big\},\qq\ae t\in[0,T].$$
Thus, $\xi(\cd)\in L^2_{\cF_T}(0,T;\dbR^m)$, and we have defined a
map $\bar\xi(\cd)\mapsto\xi(\cd)$ from $L^2_{\cF_T}(0,T;\dbR^m)$ to
itself. We now show that this map is a contraction. To this end, we
take any $\bar\xi_i(\cd)\in L^2_{\cF_T}(0,T;\dbR^m)$, $i=1,2$, and
let $\l_i(\cd\,,\cd)$ solve (\ref{4.5}). Then for almost every
$t\in[0,T]$ and any $r\in[t,T],$
$$\ba{ll}
\ns\ds\l_1(t,r)-\l_2(t,r)=\int_t^rB(t,s)^T\(\l_1(t,s)-\l_2(t,s)\)dW(s)\\
\ns\ds\qq\qq\qq\qq\q+\int_0^tA(s,t)^T\dbE_t[\bar\xi_1(s)-\bar\xi_2(s)]ds
+\int_0^t\dbE_s\[C(s,t)^T\(\bar\xi_1(s)-\bar\xi_2(s)\)\]dW(s).\ea$$
Hence,
$$\ba{ll}
\ns\ds\dbE|\l_1(t,r)-\l_2(t,r)|^2\le
K\dbE\(\int_t^r|\l_1(t,s)-\l_2(t,s)|^2ds+\int_0^t|\bar\xi_1(s)-\bar\xi_2(s)|^2ds\).\ea$$
By Grownall's inequality, we obtain
$$\dbE|\l_1(t,r)-\l_2(t,r)|^2\le K\dbE\int_0^t|\bar\xi_1(s)-\bar\xi_2(s)|^2ds,
\qq \forall r\in[t,T].$$
With the definition
$$\xi_i(t)=\l_i(t,T),\qq t\in[0,T], \ a.e.$$
one obtains
$$\dbE|\xi_1(t)-\xi_2(t)|^2\le K\dbE\int_0^t|\bar\xi_1(s)-\bar\xi_2(s)|^2ds,
\qq\ae t\in[0,T].$$
Now, for any $\m>0$, it follows that
$$\ba{ll}
\ns\ds\dbE\int_0^Te^{-\m t}|\xi_1(t)-\xi_2(t)|^2dt\le
K\dbE\int_0^Te^{-\m t}\int_0^t|\bar\xi_1(s)-\bar\xi_2(s)|^2dsdt\\
\ns\ds=K\dbE\int_0^T\(\int_s^Te^{-\m
t}dt\)|\bar\xi_1(s)-\bar\xi_2(s)|^2ds\le{K\over\m}\int_0^Te^{-\m
t}|\bar\xi_1(s)-\bar\xi_2(s)|^2ds.\ea$$
Since $K$ in the above is an absolute constant, by choosing $\m>0$
large, we get that the map $\bar\xi(\cd)\mapsto\xi(\cd)$ is a
contraction in $L^2_{\cF_T}(0,T;\dbR^m)$ with a weighted norm.
Hence, it admits a unique fixed point which is the unique solution
of (\ref{4.4}).

\ms

Next, we prove the duality relation (\ref{duality}).

\ms

For any given $\psi(\cd)\in L^2_{\cF_T}(0,T;\dbR^n)$, let
$(Y(\cd),Z(\cd\,,\cd))\in\cM^2[0,T]$ be the unique adapted
M-solution of (\ref{4.1}). Then we have
\bel{}Y(t)=\dbE Y(t)+\int_0^tZ(t,s)dW(s),\qq\ae t\in[0,T].\ee
Hence, for any $(\a(\cd),\b(\cd\,,\cd))\in
L^2_\dbF(0,T;\dbR^n)\times L^2(0,T;L^2_{\dbF}(0,T;\dbR^n))$, we
have,
\bel{4.9}\ba{ll}
\ns\ds\dbE\Big\{\int_0^T\lan Y(t),\a(t)+\int_0^t\b(t,s)dW(s)\ran
dt+\int_0^T\int_t^T\lan
Z(t,s),\b(t,s)\ran dsdt\Big\}\\
\ns\ds=\dbE\Big\{\int_0^T\lan Y(t),\a(t)\ran dt+\int_0^T\int_0^T\lan
Z(t,s),\b(t,s)\ran dsdt\Big\}.\ea\ee
On the other hand, let $\xi(\cd)\in L^2_{\cF_T}(0,T;\dbR^m)$ be the
solution to the adjoint equation (\ref{4.2}) corresponding to
$(\a(\cd),\b(\cd\,,\cd))$. Then
\bel{4.10}\ba{ll}
\ns\ds\xi(t)=\a(t)+\int_0^T\b(t,s)dW(s)+\int_t^TB(t,s)^T\dbE_s[\xi(t)]dW(s)\\
\ns\ds\qq\qq+\int_0^tA(s,t)^T\dbE_t[\xi(s)]ds+\int_0^t\dbE_s[C(s,t)^T\xi(s)]dW(s)\\
\ns\ds\qq=\a(t)+\int_0^tA(s,t)^T\dbE_t[\xi(s)]ds+\int_0^t\(\dbE_s[C(s,t)^T\xi(s)]+\b(t,s)\)dW(s)\\
\ns\ds\qq\qq+\int_t^T\(B(t,s)^T\dbE_s[\xi(t)]+\b(t,s)\)dW(s)\\
\ns\ds\qq=\dbE_t[\xi(t)]+\int_t^T\z(t,s)dW(s),\ea\ee
with
\bel{4.11}\z(t,s)=B(t,s)^T\dbE_s[\xi(t)\big]+\b(t,s),\qq
s\in[t,T],~\ae t\in[0,T].\ee
Consequently,
$$\ba{ll}
\ns\ds\dbE\int_0^T\lan\psi(t),\xi(t)\ran dt\\
\ns\ds=\dbE\Big\{\int_0^T\lan
Y(t)-\int_t^T\big[A(t,s)Y(s)+B(t,s)Z(t,s)
+C(t,s)Z(s,t)\big]ds,\xi(t)\ran dt\\
\ns\ds\qq+\int_0^T\int_t^T
\lan Z(t,s),\z(t,s)\ran dsdt\Big\}\\
\ns\ds=\dbE\Big\{\int_0^T\lan
Y(t),\xi(t)-\int_0^tA(s,t)^T\xi(s)ds\ran dt\\
\ns\ds\qq+\int_0^T\int_t^T\lan Z(t,s),\z(t,s)-B(t,s)^T\xi(t)\ran
dsdt-\int_0^T\int_0^t\lan Z(t,s),C(s,t)^T\xi(s)\ran
dsdt\Big\}\\
\ns\ds=\dbE\Big\{\int_0^T\lan
Y(t),\xi(t)-\int_0^tA(s,t)^T\xi(s)ds-\int_0^t\dbE_s\big[C(s,t)^T\xi(s)
\big]dW(s)\ran dt\\
\ns\ds\qq+\int_0^T\int_t^T\lan Z(t,s),\z(t,s)-B(t,s)^T\xi(t)\ran
dsdt\Big\}\\
\ns\ds=\dbE\Big\{\int_0^T\lan
Y(t),\int_0^t\beta(t,s)dW(s)+\alpha(t)\ran dt+\int_0^T\int_t^T\lan
Z(t,s),\beta(t,s)\ran dsdt \Big\} .\ea$$
where the last equality in the above follows from (\ref{4.10}) and
(\ref{4.11}). Thanks to (\ref{4.9}), our conclusion follows.
\endpf

\ms

Let us recall the following duality principle found in \cite{Yong
2008}, which is a corollary of Theorem 4.1.

\ms

\bf Corollary 4.2. \sl Let $X(\cd)$ be the solution to the following
FSVIE:
$$X(t)=\f(t)+\int_0^tA_0(t,s)X(s)ds+\int_0^tC_0(t,s)X(s)dW(s),\qq
t\in[0,T],$$
for uniformly bounded measurable $A_0(\cd\,,\cd)$ and
$C_0(\cd\,,\cd)$, with $s\mapsto(A_0(t,s),C_0(t,s))$ is
$\dbF$-progressively measurable. Let $(p(\cd),q(\cd\,,\cd))$ be the
adapted M-solution to the following BSVIE:
$$p(t)=\psi(t)+\int_t^T\[A_0(s,t)^Tp(s)+C_0(s,t)^Tq(s,t)\]ds-\int_t^Tq(t,s)
dW(s),\q t\in[0,T].$$
Then
\bel{}\dbE\int_0^T\lan\psi(t),X(t)\ran
dt=\dbE\int_0^T\lan\f(t),p(t)\ran dt.\ee

\rm

This is a special case of Theorem 4.1 in which
$$\a(t)=\f(t),\q\b(t,s)=0,$$
and
$$A(t,s)=A_0(s,t)^T,\q B(t,s)=0,\q C(t,s)=C_0(s,t)^T.$$
In the current case, $\xi(\cd)=X(\cd)$ is adapted. Therefore,
$$\dbE_s[C(s,t)^T\xi(s)]=\dbE_s[C_0(t,s)X(s)]=C_0(t,s)X(s).$$

\ms

\section{Proofs of Theorem 3.2 and Theorem 3.3}

\rm

This section is devoted to the proofs of Theorems 3.2 and 3.3.

\ms

As a standard step, to obtain the maximum principle we need to
obtain the variation of the state and the cost functional with
respect to the control and then use duality principle(s).

\ms

For Theorem 3.1, we have the following result.

\ms

\bf Theorem 5.1. \sl Let {\rm(H3)} and {\rm(H5)} hold and $(\bar
X(\cd),\bar Y(\cd),\bar Z(\cd\,,\cd),\bar u(\cd))$ be an optimal
4-tuple of Problem {\rm(C1)}. Then for any $v(\cd)\in\cU[0,T]$,
\bel{5.1}\ba{ll}
\ns\ds0\le\dbE\Big\{h_xX_1(T)+h_yY_1(0)+\int_0^T\int_t^T
\(f_x(t,s)X_1(s)+f_y(t,s)Y_1(s)\\
\ns\ds\qq\qq+f_z(t,s)Z_1(t,s)+f_u(t,s)[v(s)-\bar
u(s)]\)dsdt\Big\},\ea\ee
where
$$\left\{\2n\ba{ll}
\ns\ds h_x=h_x(\bar X(T),\bar Y(0)),\q h_y=h_y(\bar X(T),\bar Y(0)),\\
\ns\ds f_x(t,s)=f_x(t,s,\bar X(s),\bar Y(s),\bar Z(t,s),\bar
u(s)),\q f_y(t,s)=f_y(t,s,\bar X(s),\bar Y(s),\bar Z(t,s),\bar
u(s)),\\
\ns\ds f_z(t,s)=f_z(t,s,\bar X(s),\bar Y(s),\bar Z(t,s),\bar
u(s)),\q f_u(t,s)=f_u(t,s,\bar X(s),\bar Y(s),\bar Z(t,s),\bar
u(s)),\ea\right.$$
and $(X_1(\cd),Y_1(\cd),Z_1(\cd\,,\cd))$ is the adapted solution to
the following linear FBSVIE:
\bel{variational}\left\{\2n\ba{ll}
\ns\ds X_1(t)=\int_0^t\Big\{b_x(t,s)X_1(s)+
b_u(t,s)\big[v(s)-\bar u(s)\big]\Big\}ds\\
\ns\ds\qq\qq
+\int_0^t\Big\{\si_x(t,s)X_1(s)+\si_u(t,s)\big[v(s)-\bar u(s)\big]\Big\}dW(s),\\
\ns\ds
Y_1(t)=\psi_{x'}(t)X_1(t)+\psi_x(t)X_1(T)+\int_t^T\Big\{g_{x'}(t,s)X_1(t)
+g_x(t,s)X_1(s)+g_y(t,s)Y_1(s)\\
\ns\ds\qq\qq+g_z(t,s)Z_1(t,s)\1n+\1n g_u(t,s)\big[v(s)-\bar
u(s)\big] \Big\}ds\1n-\2n\int_t^T\2n Z_1(t,s)dW(s).\ea\right.\ee

\ms

\it Proof. \rm Let $(\bar X(\cd),\bar Y(\cd),\bar Z(\cd\,,\cd),\bar
u(\cd))$ be an optimal 4-tuple of Problem (C1). Fix any
$v(\cd)\in\cU[0,T]$ and any $\e\in(0,1)$, let
$(X^\e(\cd),Y^\e(\cd),Z^\e(\cd\,,\cd))$ be the adapted M-solution of
(\ref{FBSVIE3.1}) corresponding to
$$u^\e(\cd)\equiv\bar u(\cd)+\e[v(\cd)-\bar u(\cd)]=(1-\e)\bar u(\cd)+\e v(\cd)\in\cU[0,T].$$
Under (H3), one has
\bel{5.0}\ba{ll}
\ns\ds \sup\limits_{t\in[0,T]}\dbE|\bar
X(t)|^2+\sup\limits_{t\in[0,T]}\dbE|X^{\e}(t)|^2<\infty,\qq\forall\e\in(0,1].
\ea\ee
Note that
$$\ba{ll}
\ns\ds X^\e(t)=\f(t)+\int_0^tb(t,s,X^\e(s),u^\e(s))ds+\int_0^t\si(t,s,X^\e(s),u^\e(s))dW(s),\\
\ns\ds\bar X(t)=\f(t)+\int_0^tb(t,s,\bar X(s),\bar
u(s))ds+\int_0^t\si(t,s,\bar X(s),\bar u(s))dW(s),\ea\q\ae
t\in[0,T].$$
Thus, by Proposition 3.1, we have
\bel{}\ba{ll}
\ns\ds\sup_{t\in[0,T]}\dbE|X^\e(t)\1n-\1n \bar
X(t)|^2\1n+\2n\sup_{t\in[0,T]}\dbE|Y^\e(t)\1n-\1n\bar Y(t)|^2
\1n+\2n\sup_{t\in[0,T]}\dbE\2n\int_t^T
\2n|Z^\e(t,s)\1n-\1n\bar Z(t,s)|^2ds \\
\ns\ds\le\1n K\Big\{\2n\sup_{t\in[0,T]}\dbE\(\int_0^T\2n|
b(t,s,\bar X(s),u^\e(s))\1n-\1n b(t,s,\bar X(s),\bar u(s))|ds\)^2\\
\ns\ds\qq+\sup_{t\in[0,T]}\dbE\int_0^T|\si(t,s,\bar X(s),u^\e(s))-\si(t,s,\bar X(s),\bar u(s))|^2ds\\
\ns\ds\qq+\sup_{t\in[0,T]}\dbE\(\int_0^T|g(t,s,\bar X(t),\bar X(s),\bar Y(s),\bar Z(t,s),u^\e(s))\\
\ns\ds\qq\qq\qq\qq-g(t,s,\bar X(t),\bar X(s),\bar Y(s),\bar
Z(t,s),\bar u(s))|ds\)^2\Big\}\\
\ns\ds\le\1n K\e^2\dbE\int_0^T\2n| v(s))\1n-\1n \bar
u(s)|^2ds.\ea\ee
For (\ref{variational}), we firstly have
\bel{}\sup_{t\in[0,T]}\dbE|X_1(t)|^2\le K\dbE\int_0^T|v(t)-\bar
u(t)|^2dt.\ee
Then, by Theorem 2.3, one has
\bel{}\sup_{t\in[0,T]}\dbE|Y_1(t)|^2+\sup_{t\in[0,T]}\dbE\int_t^T|Z_1(t,s)|^2ds\le
K\dbE\int_0^T|v(t) -\bar u(t)|^2dt.\ee
Now, for $t,s\in [0,T]$, let
\bel{5.5}X_1^\e(t)={X^\e(t)-\bar X(t)\over\e},\q
Y^\e_1(t)={Y^\e(t)-\bar Y(t)\over\e},\q Z_1^\e(t,s)={Z^\e(t,s)-\bar
Z(t,s)\over\e}.\ee
Then we look at the following:
$$\ba{ll}
\ns\ds
X^\e_1(t)-X_1(t)=\int_0^t\Big\{{b(t,s,X^\e(s),u^\e(s))-b(t,s,\bar
X(s),\bar u(s))\over\e}\\
\ns\ds\qq\qq\qq\qq\qq\qq-b_x(t,s)X_1(s)-b_u(t,s)[v(s)-\bar u(s)]\Big\}ds\\
\ns\ds\qq\qq\qq\qq+\int_0^t\Big\{{\si(t,s,X^\e(s),u^\e(s))-\si(t,s,\bar
X(s),\bar u(s))\over\e}\\
\ns\ds\qq\qq\qq\qq\qq\qq-\si_x(t,s)X_1(s)-\si_u(t,s)[v(s)-\bar
u(s)]\Big\}dW(s)\\
\ns\ds=\int_0^t\Big\{b_x^\e(t,s)X^\e_1(s)-b_x(t,s)X_1(s)+\big[b_u^\e(t,s)-b_u(t,s)\big][v(s)-\bar
u(s)]\Big\}ds\\
\ns\ds\q+\int_0^t\Big\{\si_x^\e(t,s)X^\e_1(s)-\si_x(t,s)X_1(s)+\big[\si_u^\e(t,s)-\si_u(t,s)\big][v(s)-\bar
u(s)]\Big\}dW(s),\ea$$
where
$$b_x^\e(t,s)=\int_0^1b_x\big(t,s,\bar X(s)+\th[X^\e(s)-\bar X(s)],\bar
u(s)+\th\e[v(s)-\bar u(s)]\big)d\th,$$
and $b^\e_u(t,s)$, $\si^\e_x(t,s)$, $\si^\e_u(t,s)$ are defined
similarly. Then we have
$$\ba{ll}
\ns\ds\dbE|X_1^\e(t)-X_1(t)|^2\le
K\Big\{\dbE\[\int_0^t\(|b_x^\e(t,s)-b_x(t,s)|\,|X_1(s)|+|b_u^\e(t,s)-b_u(t,s)|\,|v(s)-\bar
u(s)|\)ds\]^2\\
\ns\ds\qq\qq\qq+\dbE\int_0^t\(|\si_x^\e(t,s)-\si_x(t,s)|^2|X_1(s)|^2+|\si_u^\e(t,s)-\si_u(t,s)|^2|v(s)-\bar
u(s)|^2\)ds\Big\}.\ea$$
By the dominated convergence theorem, we have
$$\left\{\ba{ll}
\ns\ds\lim_{\e\to0}\dbE|X_1^\e(t)-X_1(t)|^2=0,\qq\forall t\in[0,T],\\
\ns\ds \lim_{\e\to0}\dbE\int_0^T|X_1^\e(t)-X_1(t)|^2dt=0.
\ea\right.$$
Next,
$$\ba{ll}
\ns\ds Y^\e_1(t)-Y_1(t)={\psi(t,X^\e(t),X^\e(T))-\psi(t,\bar
X(t),\bar
X(T))\over\e}-\psi_{x'}(t)X_1(t)-\psi_x(t)X_1(T)\\
\ns\ds\q+\int_t^T\Big\{{g(t,s,X^\e(t),X^\e(s),Y^\e(s),Z^\e(t,s),u^\e(s))-g(t,s,\bar
X(t),\bar X(s),\bar Y(s),\bar Z(t,s),\bar u(s))\over\e}\\
\ns\ds\q-g_{x'}(t,s)X_1(t)-g_x(t,s)X_1(s)-g_y(t,s)Y_1(s)-g_z(t,s)Z_1(t,s)-g_u(t,s)[v(s)-\bar
u(s)]\Big\}ds\\
\ns\ds\q-\int_t^T\(Z_1^\e(t,s)-Z_1(t,s)\)dW(s)\\
\ns\ds=\psi^\e_{x'}(t)X^\e(t)-\psi_{x'}(t)X_1(t)+\psi^\e_x(t)X_1^\e(T)-\psi_x(t)X_1(T)\\
\ns\ds\q+\int_t^T\Big\{g_{x'}^\e(t,s)X_1^\e(t)-g_{x'}(t,s)X_1(t)+g_x^\e(t,s)X_1^\e(s)-g_x(t,s)X_1(s)\\
\ns\ds\q+g_y^\e(t,s)Y_1^\e(s)-g_y(t,s)Y_1(s)+g_z^\e(t,s)Z_1^\e(t,s)-g_z(t,s)Z_1(t,s)\\
\ns\ds\q+\big[g_u^\e(t,s)-g_u(t,s)\big][v(s)-\bar
u(s)]\Big\}ds-\int_t^T\(Z_1^\e(t,s)-Z_1(t,s)\)dW(s),\ea$$
where
$$\psi_{x'}^\e(t)=\int_0^1\psi_{x'}\big(t,\bar X(t)+\th[X^\e(t)-\bar X(t)],\bar X(T)+\th[X^\e(T)-\bar X(T)]\big)
d\th,$$
and $\psi_x^\e(t)$, $g_{x'}^\e(t,s)$, $g_x^\e(t,s)$, $g^\e_y(t,s)$,
$g^\e_z(t,s)$ and $g^\e_u(t,s)$ are defined similarly. Then
$$\ba{ll}
\ns\ds\dbE|Y^\e_1(t)-Y_1(t)|^2+\dbE\int_t^T|Z_1^\e(t,s)-Z_1(t,s)|^2ds\\
\ns\ds\le
K\dbE\,\Big|\psi^\e_{x'}(t)X^\e(t)-\psi_{x'}(t)X_1(t)+\psi^\e_x(t)X_1^\e(T)-\psi_x(t)X_1(T)\\
\ns\ds\qq+\int_t^T\Big\{g_{x'}^\e(t,s)X_1^\e(t)-g_{x'}(t,s)X_1(t)+g_x^\e(t,s)X_1^\e(s)-g_x(t,s)X_1(s)\\
\ns\ds\qq+\(g_y^\e(t,s)-g_y(t,s)\)Y_1(s)+\(g_z^\e(t,s)-g_z(t,s)\)Z_1(t,s)\\
\ns\ds\qq+\big[g_u^\e(t,s)-g_u(t,s)\big][v(s)-\bar
u(s)]\Big\}ds\Big|^2\ea$$
$$\ba{ll}
\ns\ds\le
K\dbE\Big\{|\psi^\e_{x'}(t)-\psi_{x'}(t)|^2|X_1(t)|^2+|\psi^\e_x(t)-\psi_x(t)|^2|X_1(T)|^2+|X_1^\e(t)-X_1(t)|^2\\
\ns\ds\qq+|X_1^\e(T)-X_1(T)|^2+\[\int_t^T\(|g_{x'}^\e(t,s)-g_{x'}(t,s)||X_1(t)|+|X_1^\e(t)-X_1(t)|\\
\ns\ds\qq+|g_x^\e(t,s)-g_x(t,s)||X_1(s)|+|X_1^\e(s)-X_1(s)|\\
\ns\ds\qq+|g_y^\e(t,s)-g_y(t,s)||Y_1(s)|+|g_z^\e(t,s)-g_z(t,s)||Z_1(t,s)|\)ds\]^2\Big\}.\ea$$
Then by the dominated convergence theorem, we have
$$\left\{\ba{ll}
\ns\ds\lim_{\e\to0}\(\dbE|Y^\e_1(t)-Y_1(t)|^2+\dbE\int_t^T|Z_1^\e(t,s)-Z_1(t,s)|^2ds\)=0,\qq\forall
t\in[0,T],\\
\ns\ds
\lim_{\e\to0}\(\dbE\int_0^T|Y^\e_1(t)-Y_1(t)|^2dt+\dbE\int_0^T\int_t^T|Z_1^\e(t,s)-Z_1(t,s)|^2dsdt\)=0
\ea\right.$$
Now, by the optimality of $(\bar X(\cd),\bar Y(\cd),\bar
Z(\cd\,,\cd),\bar u(\cd))$, we have
\bel{5.8.0}\ba{ll}
\ns\ds0\le{J_1(u^\e(\cd))-J_1(\bar u(\cd))\over\e}=\dbE\[{h(X^\e(T),Y^\e(0))-h(\bar X(T),\bar Y(0))\over\e}\\
\ns\ds\qq+\int_0^T\2n\int_t^T{f(t,s,X^\e(s),Y^\e(s),Z^\e(t,s),u^\e(s))-f(t,s,\bar
X(s),\bar Y(s),\bar Z(t,s),\bar u(s))\over\e}dsdt\]\\
\ns\ds\q=\dbE\Big\{h^\e_xX^\e_1(T)+h^\e_yY^\e_1(0)
+\int_0^T\int_t^T\(f_x^\e(t,s)X^\e_1(s)+f_y^\e(t,s)Y^\e_1(s)\\
\ns\ds\qq\qq+f_z^\e(t,s)Z^\e_1(t,s)+f_u^\e(t,s)[v(s)-\bar
u(s)]\)dsdt\Big\},\ea\ee
where
$$h^\e_x=\int_0^1h_x(\bar X(T)+\th[X^\e(T)-\bar X(T)],\bar
Y(0)+\th[Y^\e(0)-\bar Y(0)])d\th,$$
and $h^\e_y$, $f_x^\e(t,s)$, $f_y^\e(t,s)$, $f^\e_z(t,s)$, and
$f^\e_u(t,s)$ are defined similarly. Passing to the limit in
(\ref{5.8.0}), we obtain (\ref{5.1}). \endpf

\ms

We now present a lemma which will play an interesting role below.

\ms

\bf Lemma 5.2. \sl Suppose
\bel{eta_0}\eta_0=\xi+\int_0^Tg^0(s,\z_0(s))ds-\int_0^T\z_0(s)dW(s),\ee
where $\xi\in L^2_{\cF_T}(\O;\dbR^m)$, $\eta_0\in\dbR^m$ is
deterministic, $\z_0(\cd)\in L^2_\dbF(0,T;\dbR^m)$,
$g^0:[0,T]\times\dbR^m\times\O\to\dbR^m$ such that $s\mapsto
g^0(s,\z)$ is progressively measurable for all $\z\in\dbR^m$, almost
surely, and $\z\mapsto g^0(s,\z)$ is uniformly Lipschitz continuous.
Let $(\eta(\cd),\z(\cd))$ be the unique adapted solution to the
following BSDE:
\bel{BSDE5.9}\eta(t)=\xi+\int_t^Tg^0(s,\z(s))ds-\int_t^T\z(s)dW(s),\qq
t\in[0,T].\ee
Then
\bel{eta=eta}\eta_0=\eta(0),\qq\z_0(s)=\z(s),\q s\in[0,T],~\as\ee

\ms

\it Proof. \rm Under our condition, BSDE (\ref{BSDE5.9}) admits a
unique adapted solution $(\eta(\cd),\z(\cd))$. Let
$$\ba{ll}
\ns\ds\bar\eta(t)=\eta_0-\int_0^tg^0(s,\z_0(s))ds+\int_0^t\z_0(s)dW(s)\\
\ns\ds\qq=\xi+\int_t^Tg^0(s,\z_0(s))ds-\int_t^T\z_0(s)dW(s),\qq
t\in[0,T].\ea$$
Therefore, $(\bar\eta(\cd),\z_0(\cd))$ is also an adapted solution
of BSDE (\ref{BSDE5.9}). Hence, by uniqueness, one must have
$$(\bar\eta(\cd),\z_0(\cd))=(\eta(\cd),\z(\cd)).$$
Consequently, we have (\ref{eta=eta}). \endpf

\ms

We now carry out a proof of Theorem 3.2.

\ms

\it Proof of Theorem 3.2. \rm We begin with the variational
inequality (\ref{5.1}) which is rewrite here for convenience:
$$\ba{ll}
\ns\ds0\le\dbE\Big\{h_xX_1(T)+h_yY_1(0)+\int_0^T\int_t^T
\(f_x(t,s)X_1(s)+f_y(t,s)Y_1(s)\\
\ns\ds\qq\qq+f_z(t,s)Z_1(t,s)+f_u(t,s)[v(s)-\bar
u(s)]\)dsdt\Big\},\ea$$
where $(X_1(\cd),Y_1(\cd),Z_1(\cd\,,\cd))$ is the adapted solution
to the following linear FBSVIE:
$$\left\{\2n\ba{ll}
\ns\ds X_1(t)=\int_0^t\Big\{b_x(t,s)X_1(s)+
b_u(t,s)\big[v(s)-\bar u(s)\big]\Big\}ds\\
\ns\ds\qq\qq
+\int_0^t\Big\{\si_x(t,s)X_1(s)+\si_u(t,s)\big[v(s)-\bar u(s)\big]\Big\}dW(s),\\
\ns\ds
Y_1(t)=\psi_{x'}(t)X_1(t)+\psi_x(t)X_1(T)+\int_t^T\Big\{g_{x'}(t,s)X_1(t)
+g_x(t,s)X_1(s)+g_y(t,s)Y_1(s)\\
\ns\ds\qq\qq+g_z(t,s)Z_1(t,s)\1n+\1n g_u(t,s)\big[v(s)-\bar
u(s)\big] \Big\}ds\1n-\2n\int_t^T\2n Z_1(t,s)dW(s).\ea\right.$$
To obtain the maximum principle, we will use duality principles to
eliminate $Y_1(0)$, $(Y_1(\cd),Z_1(\cd\,,\cd))$, $X_1(T)$, and
$X_1(\cd)$ one by one. Therefore, there are four steps.

\ms

\it Step 1. \rm Eliminate $Y_1(0)$. We note that
$$\ba{ll}
\ns\ds Y_1(0)=\psi_x(0)X_1(T)+\int_0^T\Big\{
g_x(0,s)X_1(s)+g_y(0,s)Y_1(s)\\
\ns\ds\qq\qq+g_z(0,s)Z_1(0,s)\1n+\1n g_u(0,s)\big[v(s)-\bar
u(s)\big] \Big\}ds\1n-\2n\int_0^T\2n Z_1(0,s)dW(s).\ea$$
By Lemma 5.2, we let $(\eta_1(\cd),\z_1(\cd))$ be the adapted
solution to the following BSDE
$$\ba{ll}
\ns\ds\eta_1(t)=\psi_x(0)X_1(T)+\int_t^T\Big\{
g_x(0,s)X_1(s)+g_y(0,s)Y_1(s)\\
\ns\ds\qq\qq+g_z(0,s)\z_1(s)\1n+\1n g_u(0,s)\big[v(s)-\bar u(s)\big]
\Big\}ds\1n-\2n\int_t^T\2n\z_1(s)dW(s).\ea$$
Then one has
$$Y_1(0)=\eta_1(0),\qq Z_1(0,s)=\z_1(s),\q s\in[0,T].$$
Now, let
\bel{}\left\{\2n\ba{ll}
\ns\ds d\l(t)=g_z(0,t)^T\l(t)dW(t),\qq
t\in[0,T],\\
\ns\ds\l(0)=\dbE(h_y)^T.\ea\right.\ee
Then
$$\ba{ll}
\ns\ds\dbE\big[h_yY_1(0)\big]=\dbE\lan\l(0),\eta_1(0)\ran\\
\ns\ds=\dbE\Big\{\lan\l(T),\eta_1(T)\ran+\int_0^T\(\lan\l(s),g_x(0,s)X_1(s)
+g_y(0,s)Y_1(s)\\
\ns\ds\qq\qq\qq\qq\qq+g_z(0,s)\z_1(s)+g_u(0,s)[v(s)-\bar
u(s)]\ran-\lan g_z(0,s)^T\l(s),\z_1(s)\ran\)ds\Big\}\\
\ns\ds=\dbE\Big\{\2n\lan\l(T),\psi_x(0)X_1(T)\ran\1n+\1n\int_0^T\2n\2n\lan\l(s),
g_x(0,s)X_1(s)\1n+\1n g_y(0,s)Y_1(s)\1n+\1n g_u(0,s)[v(s)\1n-\1n\bar
u(s)]\ran ds\Big\}.\ea$$
Hence,
$$\ba{ll}
\ns\ds J_1(\bar u(\cd);v(\cd))=\dbE\Big\{h^T_xX_1(T)+h^T_yY_1(0)
+\int_0^T\int_t^T\(f_x(t,s)^TX_1(s)+f_y(t,s)^TY_1(s)\\
\ns\ds\qq\qq\qq\qq+f_z(t,s)^TZ_1(t,s)+f_u(t,s)^T[v(s)-\bar
u(s)]\)dsdt\Big\}\\
\ns\ds\q=\dbE\Big\{\lan h_x^T+\psi_x(0)^T\l(T),X_1(T)\ran
+\int_0^T\(\lan g_x(0,s)^T\l(s)+\int_0^sf_x(t,s)^Tdt,X_1(s)\ran\\
\ns\ds\qq\q+\lan g_u(0,s)^T\l(s)+\int_0^sf_u(t,s)^Tdt,v(s)-\bar
u(s)\ran\)ds\\
\ns\ds\qq\q+\int_0^T\2n\lan g_y(0,s)^T\1n\l(s)\1n+\2n\int_0^s\2n
f_y(t,s)^Tdt,Y_1(s)\ran ds+\2n\int_0^T\2n\int_t^T\2n\lan
f_z(t,s)^T,Z_1(t,s)\ran dsdt\Big\}.\ea$$
So, $Y_1(0)$ is eliminated.

\ms

\it Step 2. \rm Eliminate $(Y_1(\cd),Z_1(\cd\,,\cd))$. We apply
Theorem 4.1 with the following:
$$A(t,s)=g_y(t,s),\qq B(t,s)=g_z(t,s),\qq C(t,s)=0,$$
and
$$\ba{ll}
\ns\ds\psi(t)=\psi_{x'}(t)X_1(t)+\psi_x(t)X_1(T)+\int_t^T\Big\{g_{x'}(t,s)X_1(t)
+g_x(t,s)X_1(s)\\
\ns\ds\qq\qq+g_u(t,s)\big[v(s)-\bar u(s)\big]\1n\Big\}ds,\\
\ns\ds\a(t)=g_y(0,t)^T\l(t)+\int_0^tf_y(r,t)^Tdr,\qq\b(t,s)=f_z(t,s)^T.\ea$$
Therefore, let $\xi(\cd)$ solve the following:
\bel{}\ba{ll}
\ns\ds\xi(t)=g_y(0,t)^T\l(t)+\int_0^tf_y(r,t)^Tdr
+\int_t^Tf_z(t,r)^TdW(r)+\int_t^Tg_z(t,s)^T\dbE_s[\xi(t)]dW(s)\\
\ns\ds\qq\qq+\int_0^tg_y(s,t)^T\dbE_t[\xi(s)]ds, \qq\ae t\in[0,T].\
\ea\ee
Then by Theorem 4.1,
\bel{}\ba{ll}
\ns\ds\dbE\Big\{\int_0^T\lan
g_y(0,t)^T\l(t)+\int_0^tf_y(r,t)^Tdr,Y_1(t)\ran
dt+\int_0^T\int_t^T\lan f_z(t,s)^T,Z_1(t,s)\ran dsdt\Big\}\\
\ns\ds=\dbE\int_0^T\lan\(\psi_{x'}(t)+\int_t^Tg_{x'}(t,s)ds\)X_1(t)
+\psi_x(t)X_1(T)+\int_t^Tg_x(t,s)X_1(s)ds\\
\ns\ds\qq+\int_t^Tg_u(t,s)[v(s)-\bar u(s)]ds,\xi(t)\ran dt.\ea\ee
Hence,
$$\ba{ll}
\ns\ds J_1(\bar u(\cd);v(\cd))=\dbE\Big\{\lan
h_x^T+\psi_x(0)^T\l(T),X_1(T)\ran
+\int_0^T\(\lan g_x(0,s)^T\l(s)+\int_0^sf_x(t,s)^Tdt,X_1(s)\ran\\
\ns\ds\qq\q+\lan g_u(0,s)^T\l(s)+\int_0^sf_u(t,s)^Tdt,v(s)-\bar
u(s)\ran\)ds\\
\ns\ds\qq\q+\int_0^T\2n\lan g_y(0,s)^T\1n\l(s)\1n+\2n\int_0^s\2n
f_y(t,s)^T\1n dt,Y_1(s)\ran ds+ \2n\int_0^T\2n\int_t^T\2n\lan
f_z(t,s)^T,Z_1(t,s)\ran dsdt\Big\}\\
\ns\ds=\dbE\Big\{\lan h_x^T+\psi_x(0)^T\l(T),X_1(T)\ran
+\int_0^T\(\lan g_x(0,s)^T\l(s)+\int_0^sf_x(t,s)^Tdt,X_1(s)\ran\\
\ns\ds\qq\q+\lan g_u(0,s)^T\l(s)+\int_0^sf_u(t,s)^Tdt,v(s)-\bar
u(s)\ran\)ds\\
\ns\ds\qq\q+\int_0^T\lan\(\psi_{x'}(t)+\int_t^Tg_{x'}(t,s)ds\)X_1(t)
+\psi_x(t)X_1(T)+\int_t^Tg_x(t,s)X_1(s)ds\\
\ns\ds\qq\q+\int_t^Tg_u(t,s)[v(s)-\bar u(s)]ds,\xi(t)\ran
dt\Big\}\\
\ns\ds=\dbE\Big\{\lan
h_x^T+\psi_x(0)^T\l(T)+\int_0^T\psi_x(t)^T\xi(t)dt,X_1(T)\ran
+\int_0^T\[\lan
g_x(0,s)^T\l(s)+\int_0^sf_x(t,s)^Tdt\\
\ns\ds\qq\q+\(\psi_{x'}(s)^T
+\int_s^Tg_{x'}(s,t)^Tdt\)\xi(s)+\int_0^sg_x(t,s)^T\xi(t)
dt,X_1(s)\ran\\
\ns\ds\qq\q+\lan
g_u(0,s)^T\l(s)+\int_0^sf_u(t,s)^Tdt+\int_0^sg_u(t,s)^T\xi(t)dt,v(s)-\bar
u(s)\ran\]ds\Big\}.\ea$$
Thus, $(Y_1(\cd),Z_1(\cd\,,\cd))$ is eliminated.

\ms

\it Step 3. \rm Eliminate $X_1(T)$. Let $(\m(\cd),\n(\cd))$ be the
adapted solution to the following BSDE:
$$\m(t)=h_x^T+\psi_x(0)^T\l(T)+\int_0^T\psi_x(r)^T\xi(r)dr-\int_t^T\n(s)dW(s),
\qq t\in[0,T].$$
Note that
$$\ba{ll}
\ns\ds X_1(T)=\int_0^T\(b_x(T,s)X_1(s)+b_u(T,s)[v(s)-\bar u(s)]\)ds\\
\ns\ds\qq\qq+\int_0^T\(\si_x(T,s)X_1(s)+\si_u(T,s)[v(s)-\bar
u(s)]\)dW(s).\ea$$
Thus,
$$\ba{ll}
\ns\ds\dbE\Big\{\lan
h_x^T+\psi_x(0)^T\l(T)+\int_0^T\psi_x(r)^T\xi(r)dr,X_1(T)\ran\Big\}\\
\ns\ds=\dbE\Big\{\int_0^T\2n\lan
h_x^T\1n+\1n\psi_x(0)^T\l(T)\1n+\2n\int_0^T\2n\psi_x(r)^T\xi(r)dr,
b_x(T,s)X_1(s)\1n+\1n b_u(T,s)[v(s)\1n-\1n\bar u(s)]\ran ds\\
\ns\ds\q+\int_0^T\lan\n(s),\si_x(T,s)X_1(s)+\si_u(T,s)[v(s)-\bar u(s)]\ran ds\Big\}\\
\ns\ds=\dbE\int_0^T\Big\{\lan b_x(T,s)^T\(
h_x^T+\psi_x(0)^T\l(T)+\int_0^T\psi_x(r)^T\xi(r)dr\)+\si_x(T,s)^T\n(s),X_1(s)\ran\\
\ns\ds\q+\lan b_u(T,s)^T\(
h_x^T+\psi_x(0)^T\l(T)+\int_0^T\psi_x(t)^T\xi(t)dt\)\\
\ns\ds\q+\si_u(T,s)^T\n(s),v(s)\1n-\1n\bar u(s)\ran\Big\} ds.\ea$$
Consequently,
$$\ba{ll}
\ns\ds J_1(\bar u(\cd);v(\cd))=\dbE\int_0^T\Big\{\lan b_x(T,s)^T\(
h_x^T+\psi_x(0)^T\l(T)+\int_0^T\psi_x(r)^T\xi(r)dr\)\\
\ns\ds\qq\qq\qq\q+\si_x(T,s)^T\n(s)+g_x(0,s)^T\l(s)+\int_0^sf_x(r,s)^Tdr\\
\ns\ds\qq\qq\qq\q+\(\psi_{x'}(s)^T
+\int_s^Tg_{x'}(s,r)^Tdr\)\xi(s)+\int_0^sg_x(r,s)^T\xi(r)
dr,X_1(s)\ran\\
\ns\ds\qq\qq\qq\q+\lan b_u(T,s)^T\(
h_x^T+\psi_x(0)^T\l(T)+\int_0^T\psi_x(r)^T\xi(r)dr\)+\si_u(T,s)^T\n(s)\\
\ns\ds\qq\qq\qq\q+
g_u(0,s)^T\l(s)+\int_0^sf_u(r,s)^Tdr+\int_0^sg_u(r,s)^T\xi(r)dr,v(s)-\bar
u(s)\ran\Big\}ds.\ea$$
Thus, $X_1(T)$ is eliminated.

\ms

\it Step 4. \rm Eliminate $X_1(\cd)$. We apply Corollary 4.2 with

$$A_0(t,s)=b_x(t,s),\qq C_0(t,s)=\si_x(t,s),$$
and
$$\ba{ll}
\ns\ds\f(t)=\int_0^tb_u(t,s)[v(s)-\bar u(s)]ds+\int_0^t\si_u(t,s)[v(s)-\bar u(s)]dW(s),\\
\ns\ds\psi(t)=b_x(T,t)^T\(
h_x^T+\psi_x(0)^T\l(T)+\int_0^T\psi_x(r)^T\xi(r)dr\)\\
\ns\ds\qq\qq+\si_x(T,t)^T\n(t)+g_x(0,t)^T\l(t)+\int_0^tf_x(r,t)^Tdr\\
\ns\ds\qq\qq+\(\psi_{x'}(t)^T
+\int_s^Tg_{x'}(t,r)^Tdr\)\xi(t)+\int_0^tg_x(r,t)^T\xi(r) dr.\ea$$
Therefore, we let $(p(\cd),q(\cd\,,\cd))$ be the adapted M-solution
to the following BSVIE:
$$\ba{ll}
\ns\ds p(t)=b_x(T,t)^T\(
h_x^T+\psi_x(0)^T\l(T)+\int_0^T\psi_x(r)^T\xi(r)dr\)\\
\ns\ds\qq\qq+\si_x(T,t)^T\n(t)+g_x(0,t)^T\l(t)+\int_0^tf_x(r,t)^Tdr\\
\ns\ds\qq\qq+\(\psi_{x'}(t)^T
+\int_t^Tg_{x'}(t,r)^Tdr\)\xi(t)+\int_0^tg_x(r,t)^T\xi(r)dr\\
\ns\ds\qq\qq+\int_t^T\(b_u(s,t)^Tp(s)+\si_u(s,t)^Tq(s,t)\)ds-\int_t^Tq(t,s)dW(s).\ea$$
Then, by Corollary 4.2, one obtains
$$\ba{ll}
\ns\ds\dbE\int_0^T\lan b_x(T,t)^T\(
h_x^T+\psi_x(0)^T\l(T)+\int_0^T\psi_x(r)^T\xi(r)dr\)\\
\ns\ds\qq\qq+\si_x(T,t)^T\n(t)+g_x(0,t)^T\l(t)+\int_0^tf_x(r,t)^Tdr\\
\ns\ds\qq\qq+\(\psi_{x'}(t)^T
+\int_t^Tg_{x'}(t,r)^Tdr\)\xi(t)+\int_0^tg_x(r,t)^T\xi(r)dr,X_1(t)\ran
dt\\
\ns\ds=\dbE\int_0^T\lan p(t),\int_0^tb_u(t,s)[v(s)-\bar u(s)]ds
+\int_0^t\si_u(t,s)[v(s)-\bar u(s)]dW(s)\ran\\
\ns\ds=\dbE\int_0^T\int_s^T\lan
b_u(t,s)^Tp(t)+\si_u(t,s)^Tq(t,s),v(s)-\bar u(s)\ran dtds.\ea$$
Consequently,
$$\ba{ll}
\ns\ds J_1(\bar
u(\cd);v(\cd))=\dbE\int_0^T\lan\int_s^T\big[b_u(t,s)^Tp(t)+\si_u(t,s)^T
q(t,s)\big]dt\\
\ns\ds\qq\qq\qq\qq+b_u(T,s)^T\m(T)+\si_u(T,s)^T\n(s)+
g_u(0,s)^T\l(s)\\
\ns\ds\qq\qq\qq\qq+\int_0^sf_u(t,s)^Tdt+\int_0^sg_u(t,s)^T\xi(t)dt,v(s)-\bar
u(s)\ran ds.\ea$$
Then the maximum principle follows. \endpf

\ms

Now, we look at Problem (C2). For that, we have the following
result.

\ms

\bf Theorem 5.3. \sl Let {\rm(H4)} and {\rm(H6)} hold. Let $(\bar
X(\cd),\bar Y(\cd),\bar Z(\cd\,,\cd),\bar u(\cd))$ be an optimal
4-tuple of Problem {\rm(C2)}. Then for any $v(\cd)\in\cU[0,T]$,
\bel{5.10}\ba{ll}
\ns\ds0\le\dbE\Big\{h_xX_1(T)+\int_0^T\[\(\int_0^Tf_x(t,s)^Tdt\)X_1(s)+\(h^T_y+\int_0^Tf_y(t,s)^Tdt\)
Y_1(s)\\
\ns\ds\qq+\(\int_0^T\2n f_u(t,s)^Tdt\)[v(s)-\bar u(s)]\]ds
+\2n\int_0^T\2n\int_0^T\2nf_z(t,s)^TZ_1(t,s)dsdt\Big\},\ea\ee
with
$$h_x=h_x\(\bar X(T),\int_0^T\bar Y(t)dt\),\q h_y=h_y\(\bar X(T),\int_0^T\bar Y(t)dt\),$$
etc, and $(X_1(\cd),Y_1(\cd),Z_1(\cd\,,\cd))$ is the adapted
solution to the following linear FBSVIE:
\bel{variational2.0}\left\{\2n\ba{ll}
\ns\ds X_1(t)=\int_0^t\Big\{b_x(t,s)X_1(s)+
b_u(t,s)\big[v(s)-\bar u(s)\big]\Big\}ds\\
\ns\ds\qq\qq
+\int_0^t\Big\{\si_x(t,s)X_1(s)+\si_u(t,s)\big[v(s)-\bar u(s)\big]\Big\}dW(s),\\
\ns\ds
Y_1(t)=\psi_{x'}(t)X_1(t)+\psi_x(t)X_1(T)+\int_t^T\Big\{g_{x'}(t,s)X_1(t)
+g_x(t,s)X_1(s)+g_y(t,s)Y_1(s)\\
\ns\ds\qq\qq+g_z(t,s)Z_1(t,s)\1n+\1n g_{z'}(t,s)Z_1(s,t)\1n+\1n
g_u(t,s)\big[v(s)-\bar u(s)\big]\1n\Big\}ds\1n-\2n\int_t^T\2n
Z_1(t,s)dW(s),\ea\right.\ee
with for example,
\bel{5.6}\ba{ll}
\ns\ds b_x(t,s)=b_x(t,s,\bar X(s),\bar
u(s)),\q\si_x(t,s)=\si_x(t,s,\bar X(s),\bar u(s)),\\
\ns\ds\psi_{x}(t)=\psi_x(t,\bar X(t),\bar X(T)),\q
g_x(t,s)=g_x(t,s,\bar X(t),\bar X(s),\bar Y(s),\bar Z(t,s),\bar
Z(s,t),\bar u(s)). \ea \ee
and so on.

\ms

\it Proof. \rm For optimal 4-tuple $(\bar X(\cd),\bar Y(\cd),\bar
Z(\cd\,,\cd),\bar u(\cd))$ of Problem (C2), fix any
$v(\cd)\in\cU[0,T]$, for any $\e\in(0,1)$, we define
$(X^\e(\cd),Y^\e(\cd),Z^\e(\cd\,,\cd))$ similar to the proof of
Theorem 5.1. Then
\bel{}\ba{ll}
\ns\ds\sup_{t\in[0,T]}\dbE|X^\e(t)\1n-\1n\bar
X(t)|^2\1n+\1n\dbE\int_0^T|Y^\e(t)\1n-\1n\bar
Y(t)|^2dt\1n+\1n\dbE\int_0^T\2n\int_0^T
\1n|Z^\e(t,s)\1n-\1n\bar Z(t,s)|^2dsdt\\
\ns\ds\le\1n K\Big\{\2n\sup_{t\in[0,T]}\dbE\(\int_0^T\2n|
b(t,s,\bar X(s),u^\e(s))\1n-\1n b(t,s,\bar X(s),\bar u(s))|ds\)^2\\
\ns\ds\qq+\sup_{t\in[0,T]}\dbE\int_0^T|\si(t,s,\bar X(s),u^\e(s))-\si(t,s,\bar X(s),\bar u(s))|^2ds\\
\ns\ds\qq+\dbE\int_0^T\(\int_0^T|g(t,s,\bar X(t),\bar X(s),\bar Y(s),\bar Z(t,s),\bar Z(s,t),u^\e(s))\\
\ns\ds\qq\qq\qq\qq-g(t,s,\bar X(t),\bar X(s),\bar Y(s),\bar
Z(t,s),\bar Z(s,t),\bar u(s))|ds\)^2\Big\}\\
\ns\ds\le\1n K\e^2\dbE\int_0^T\2n| v(s))\1n-\1n \bar
u(s)|^2ds.\ea\ee
For $t,s\in [0,T]$, we still let
\bel{5.5}X_1^\e(t)={X^\e(t)-\bar X(t)\over\e},\q
Y^\e_1(t)={Y^\e(t)-\bar Y(t)\over\e},\q Z_1^\e(t,s)={Z^\e(t,s)-\bar
Z(t,s)\over\e}.\ee
For the variational system (\ref{variational2.0}), we have
\bel{}\sup_{t\in[0,T]}\dbE|X_1(t)|^2\le K\dbE\int_0^T|v(t)-\bar
u(t)|^2dt.\ee
Consequently, by Theorem 2.3,
\bel{}\dbE\int_0^T|Y_1(t)|^2+\dbE\int_0^T\int_0^T|Z_1(t,s)|^2dsdt\le
K\dbE\int_0^T|v(t)-\bar u(t)|^2dt.\ee
We also have
$$\left\{\ba{ll}
\ns\ds \lim_{\e\to0}\dbE|X_1^\e(t)-X_1(t)|^2=0,\qq\forall
t\in[0,T],\\
\ns\ds
\lim_{\e\to0}\dbE\int_0^T|X_1^\e(t)-X_1(t)|^2dt=0.\ea\right.$$
Next,
$$\ba{ll}
\ns\ds Y^\e_1(t)-\bar Y_1(t)={\psi(t,X^\e(t),X^\e(T))-\psi(t,\bar
X(t),\bar
X(T))\over\e}-\psi_{x'}(t)X_1(t)-\psi_x(t)X_1(T)\\
\ns\ds\q+\int_t^T\Big\{{1\over\e}\[g(t,s,X^\e(t),X^\e(s),Y^\e(s),Z^\e(t,s),Z^\e(s,t),u^\e(s))\\
\ns\ds\qq\qq\qq-g(t,s,\bar X(t),\bar X(s),\bar Y(s),\bar Z(t,s),\bar
Z(s,t),\bar
u(s))\]-g_{x'}(t,s)X_1(t)-g_x(t,s)X_1(s)\\
\ns\ds\qq\qq\qq-g_y(t,s)Y_1(s)-g_z(t,s)Z_1(t,s)-g_{z'}(t,s)Z_1(s,t)-g_u(t,s)[v(s)-\bar
u(s)]\Big\}ds\\
\ns\ds\q-\int_t^T\(Z_1^\e(t,s)-Z_1(t,s)\)dW(s)\\
\ns\ds=\psi^\e_{x'}(t)X^\e(t)-\psi_{x'}(t)X_1(t)+\psi^\e_x(t)X_1^\e(T)-\psi_x(t)X_1(T)\\
\ns\ds\q+\int_t^T\Big\{g_{x'}^\e(t,s)X_1^\e(t)-g_{x'}(t,s)X_1(t)+g_x^\e(t,s)X_1^\e(s)-g_x(t,s)X_1(s)\\
\ns\ds\q+g_y^\e(t,s)Y_1^\e(s)-g_y(t,s)Y_1(s)+g_z^\e(t,s)Z_1^\e(t,s)-g_z(t,s)Z_1(t,s)\\
\ns\ds\q+g_{z'}(t,s)Z^\e_1(s,t)-g_{z'}(t,s)Z_1(s,t)+\big[g_u^\e(t,s)-g_u(t,s)\big][v(s)-\bar
u(s)]\Big\}ds\\
\ns\ds\q-\int_t^T\(Z_1^\e(t,s)-Z_1(t,s)\)dW(s),\ea$$
where
$$\psi_{x'}^\e(t)=\int_0^1\psi_{x'}\big(t,\bar X(t)+\th[X^\e(t)-\bar X(t)],\bar X(T)+\th[X^\e(T)-\bar X(T)]\big)
d\th,$$
and $\psi_x^\e(t)$, $g_{x'}^\e(t,s)$, $g_x^\e(t,s)$, $g^\e_y(t,s)$,
$g^\e_z(t,s)$, $g^\e_{z'}(t,s)$ and $g^\e_u(t,s)$ are defined
similarly. Then
$$\ba{ll}
\ns\ds\dbE\int_0^T|Y^\e_1(t)-Y_1(t)|^2dt+\dbE\int_0^T\int_0^T|Z_1^\e(t,s)-Z_1(t,s)|^2dsdt\\
\ns\ds\le
K\dbE\int_0^T\Big|\psi^\e_{x'}(t)X^\e(t)-\psi_{x'}(t)X_1(t)+\psi^\e_x(t)X_1^\e(T)-\psi_x(t)X_1(T)\\
\ns\ds\qq+\int_t^T\[g_{x'}^\e(t,s)X_1^\e(t)-g_{x'}(t,s)X_1(t)+g_x^\e(t,s)X_1^\e(s)-g_x(t,s)X_1(s)\\
\ns\ds\qq+\(g_y^\e(t,s)-g_y(t,s)\)Y_1(s)+\(g_z^\e(t,s)-g_z(t,s)\)Z_1(t,s)\\
\ns\ds\qq+\(g^\e_{z'}(t,s)-g_{z'}(t,s)\)Z_1(s,t)+\big[g_u^\e(t,s)-g_u(t,s)\big][v(s)-\bar
u(s)]\]ds\Big|^2dt\\
\ns\ds\le
K\dbE\int_0^T\Big\{|\psi^\e_{x'}(t)-\psi_{x'}(t)|^2|X_1(t)|^2+|\psi^\e_x(t)-\psi_x(t)|^2|X_1(T)|^2+|X_1^\e(t)-X_1(t)|^2\\
\ns\ds\qq+|X_1^\e(T)-X_1(T)|^2+\[\int_t^T\(|g_{x'}^\e(t,s)-g_{x'}(t,s)||X_1(t)|+|X_1^\e(t)-X_1(t)|\\
\ns\ds\qq+|g_x^\e(t,s)-g_x(t,s)||X_1(s)|+|X_1^\e(s)-X_1(s)|+|g_y^\e(t,s)-g_y(t,s)||Y_1(s)|\\
\ns\ds\qq+|g_z^\e(t,s)-g_z(t,s)||Z_1(t,s)|+|g_{z'}^\e(t,s)-g_{z'}(t,s)||Z_1(s,t)|\)ds\]^2\Big\}dt.\ea$$
Then by the dominated convergence theorem, we have
$$\lim_{\e\to0}\(\dbE\int_0^T|Y^\e_1(t)-Y_1(t)|^2dt+\dbE\int_0^T\int_0^T|Z_1^\e(t,s)-Z_1(t,s)|^2dsdt\)=0.$$
Now, by the optimality of $(\bar X(\cd),\bar Y(\cd),\bar
Z(\cd\,,\cd),\bar u(\cd))$, we have
\bel{5.8}\ba{ll}
\ns\ds0\le{J_2(u^\e(\cd))-J_2(\bar
u(\cd))\over\e}={1\over\e}\dbE\[h\(X^\e(T),\int_0^TY^\e(t)dt\)-h\(\bar
X(T),
\int_0^T\bar Y(t)dt\)\\
\ns\ds\qq+\int_0^T\2n\int_0^T\(f(t,s,X^\e(s),Y^\e(s),Z^\e(t,s),u^\e(s))\\
\ns\ds\qq\qq-f(t,s,\bar
X(s),\bar Y(s),\bar Z(t,s),\bar u(s))\)dsdt\]\\
\ns\ds\q=\dbE\Big\{h^\e_xX^\e_1(T)+h^\e_y\int_0^TY^\e_1(t)dt+\int_0^T\int_0^T\(f_x^\e(t,s)X^\e_1(s)+f_y^\e(t,s)Y^\e_1(s)\\
\ns\ds\qq+f_z^\e(t,s)Z^\e_1(t,s)+f_u^\e(t,s)[v(s)-\bar
u(s)]\)ds\Big\},\ea\ee
where
$$h^\e_x=\int_0^1h_x\(\bar X(T)+\th[X^\e(T)-\bar X(T)],\int_0^T\big\{\bar
Y(t)+\th[Y^\e(t)-\bar Y(t)]\big\}dt\)d\th,$$
and $h^\e_y$, $f_x^\e(t,s)$, $f_y^\e(t,s)$, $f^\e_z(t,s)$, and
$f^\e_u(t,s)$ are defined similarly. Passing to the limit in
(\ref{5.8}), by the dominated convergence theorem, we have
(\ref{5.10}). \endpf

\ms

Now let us present a proof of Theorem 3.3.

\ms

\it Proof of Theorem 3.3. \rm Similar to the proof of Theorem 3.2,
we begin with (\ref{5.10}) which is written here
$$\ba{ll}
\ns\ds0\le\dbE\Big\{h_xX_1(T)+\int_0^T\[\(\int_0^Tf_x(t,s)dt\)X_1(s)
+\(h_y+\int_0^Tf_y(t,s)dt\)
Y_1(s)\\
\ns\ds\qq+\(\int_0^T\2n f_u(t,s)dt\)[v(s)-\bar u(s)]\]ds
+\2n\int_0^T\2n\int_0^T\2nf_z(t,s)Z_1(t,s)dsdt\Big\}\equiv J_2(\bar
u(\cd);v(\cd)),\ea$$
with $(X_1(\cd),Y_1(\cd),Z_1(\cd\,,\cd))$ being the adapted
M-solution to the following linear FBSVIE:
\bel{variational2}\left\{\2n\ba{ll}
\ns\ds X_1(t)=\int_0^t\Big\{b_x(t,s)X_1(s)+
b_u(t,s)\big[v(s)-\bar u(s)\big]\Big\}ds\\
\ns\ds\qq\qq
+\int_0^t\Big\{\si_x(t,s)X_1(s)+\si_u(t,s)\big[v(s)-\bar u(s)\big]\Big\}dW(s),\\
\ns\ds
Y_1(t)=\psi_{x'}(t)X_1(t)+\psi_x(t)X_1(T)+\int_t^T\Big\{g_{x'}(t,s)X_1(t)
+g_x(t,s)X_1(s)+g_y(t,s)Y_1(s)\\
\ns\ds\qq\qq+g_z(t,s)Z_1(t,s)\1n+\1n g_{z'}(t,s)Z_1(s,t)\1n+\1n
g_u(t,s)\big[v(s)-\bar u(s)\big]\1n\Big\}ds\1n-\2n\int_t^T\2n
Z_1(t,s)dW(s).\ea\right.\ee
In the current case, we are going to eliminate
$(Y_1(\cd),Z_1(\cd\,,\cd))$, $X_1(T)$ and $X_1(\cd)$ consecutively.
Therefore, we will have three steps.

\ms

\it Step 1. \rm Eliminate $(Y_1(\cd),Z_1(\cd\,,\cd))$. To this end,
we apply Theorem 4.1 with
$$A(t,s)=g_y(t,s),\qq B(t,s)=g_z(t,s),\qq C(t,s)=g_{z'}(t,s),$$
and
$$\ba{ll}
\ns\ds\psi(t)=\psi_{x'}(t)X_1(t)+\psi_x(t)X_1(T)+\int_t^T\Big\{g_{x'}(t,s)X_1(t)
+g_x(t,s)X_1(s)\\
\ns\ds\qq\qq+g_u(t,s)\big[v(s)-\bar u(s)\big]\1n\Big\}ds,\\
\ns\ds\a(t)=h_y^T+\int_0^Tf_y(t,s)^Tdt,\qq\b(t,s)=f_z(t,s)^T,\ea$$
Therefore, we let $\xi(\cd)$ solve the following:
\bel{}\ba{ll}
\ns\ds\xi(t)=h_y^T+\int_0^Tf_y(r,t)^Tdr+\int_0^T f_z(t,s)^TdW(s)
+\int_t^Tg_z(t,s)^T\dbE_s[\xi(t)]dW(s)\\
\ns\ds\qq\qq+\int_0^tg_y(s,t)^T\dbE_t[\xi(s)]ds+\int_0^t\dbE_s[g_{z'}(s,t)^T\xi(s)]dW(s),
\qq\ae t\in[0,T].\ea\ee
Then by Theorem 4.1, one has
$$\ba{ll}
\ns\ds\dbE\Big\{\int_0^T\[\(h_y+\int_0^Tf_y(t,s)dt\)Y_1(s)ds
+\2n\int_0^T\2n\int_0^T\2nf_z(t,s)
Z_1(t,s)dsdt\Big\}\\
\ns\ds\equiv\dbE\Big\{\int_0^T\lan\a(s),Y_1(s)\ran
ds+\int_0^T\int_0^T\lan\b(t,s),Z_1(t,s)\ran
dsdt=\dbE\int_0^T\lan\psi(t),\xi(t)\ran dt\\
\ns\ds=\dbE\int_0^T\lan\psi_{x'}(t)X_1(t)+\psi_x(t)X_1(T)+\int_t^T\(g_{x'}(t,s)X_1(t)
+g_x(t,s)X_1(s)\\
\ns\ds\qq\qq+g_u(t,s)\big[v(s)-\bar u(s)\big]\)ds,\xi(t)\ran
dt\\
\ns\ds=\dbE\Big\{\lan\int_0^T\psi_x(t)^T\xi(t)dt,X_1(T)\ran\1n+\1n\int_0^T\2n\lan\(\psi_{x'}(t)\1n
+\1n\int_t^T\2n g_{x'}(t,s)ds\)^T\2n\xi(t)\1n+\2n\int_0^t\2n g_x(s,t)^T\xi(s)ds,X_1(t)\ran dt\\
\ns\ds\qq\qq+\int_0^T\lan\int_0^tg_u(s,t)^T\xi(s)ds,v(t)-\bar
u(t)\ran dt\Big\},\ea$$
which leads to
$$\ba{ll}
\ns\ds J_2(\bar
u(\cd);v(\cd))=\dbE\Big\{h_xX_1(T)+\int_0^T\[\(\int_0^Tf_x(t,s)dt\)
X_1(s)+\(h_y+\int_0^Tf_y(t,s)dt\)Y_1(s)\\
\ns\ds\qq+\(\int_0^T\2n f_u(t,s)dt\)[v(s)-\bar u(s)]\]ds
+\2n\int_0^T\2n\int_0^T\2n f_z(t,s)Z_1(t,s)dsdt\Big\}\\
\ns\ds=\dbE\Big\{h_xX_1(T)+\int_0^T\[\(\int_0^Tf_x(t,s)dt\)X_1(s)
+\(\int_0^T\2n f_u(t,s)dt\)[v(s)-\bar u(s)]\]ds\\
\ns\ds\qq+\lan\int_0^T\psi_x(t)^T\xi(t)dt,X_1(T)\ran\1n+\1n\int_0^T\2n\lan\(\psi_{x'}(t)\1n
+\1n\int_t^T\2n g_{x'}(t,s)ds\)^T\2n\xi(t)\1n+\2n\int_0^t\2n g_x(s,t)^T\xi(s)ds,X_1(t)\ran dt\\
\ns\ds\qq+\int_0^T\lan\int_0^tg_u(s,t)^T\xi(s)ds,v(t)-\bar
u(t)\ran dt\Big\}\\
\ns\ds=\dbE\Big\{\lan
h^T_x+\int_0^T\psi_x(r)^T\xi(r)dr,X_1(T)\ran\\
\ns\ds\qq+\int_0^T\lan\int_0^Tf_x(r,t)^Tdr+\(\psi_{x'}(t)\1n
+\1n\int_t^T\2n g_{x'}(t,s)ds\)^T\2n\xi(t)\1n+\2n\int_0^t\2n g_x(s,t)^T\xi(s)ds,X_1(t)\ran dt\\
\ns\ds\qq+\int_0^T\lan\int_0^tg_u(s,t)^T\xi(s)ds+\int_0^Tf_u(r,t)^Tdr,v(t)-\bar
u(t)\ran dt\Big\}.\ea$$
This finishes the elimination of $(Y_1(\cd),Z_1(\cd\,,\cd))$.

\ms

\it Step 2. \rm Elimination of $X_1(T)$. We let $(\m(\cd),\n(\cd))$
be the adapted solution to the following BSDE:
$$\m(t)=h_x^T+\int_0^T\psi_x(r)^T\xi(r)dr-\int_t^T\n(s)dW(s),
\qq t\in[0,T].$$
Note that
$$\ba{ll}
\ns\ds X_1(T)=\int_0^T\(b_x(T,s)X_1(s)+b_u(T,s)[v(s)-\bar u(s)]\)ds\\
\ns\ds\qq\qq+\int_0^T\(\si_x(T,s)X_1(s)+\si_u(T,s)[v(s)-\bar
u(s)]\)dW(s).\ea$$
Thus,
$$\ba{ll}
\ns\ds\dbE\Big\{\lan
h_x^T+\int_0^T\psi_x(r)^T\xi(r)dr,X_1(T)\ran\Big\}=\dbE\lan\m(T),X_1(T)\ran\\
\ns\ds=\dbE\Big\{\int_0^T\2n\lan\m(T),
b_x(T,s)X_1(s)\1n+\1n b_u(T,s)[v(s)\1n-\1n\bar u(s)]\ran ds\\
\ns\ds\q+\int_0^T\lan\n(s),\si_x(T,s)X_1(s)+\si_u(T,s)[v(s)-\bar u(s)]\ran ds\Big\}\\
\ns\ds=\dbE\Big\{\int_0^T\lan b_x(T,s)^T\m(T)+\si_x(T,s)^T\n(s),X_1(s)\ran\\
\ns\ds\qq\q+\lan b_u(T,s)^T\m(T)+\si_u(T,s)^T\n(s),v(s)\1n-\1n\bar
u(s)\ran ds\Big\}.\ea$$
Consequently,
$$\ba{ll}
\ns\ds J_2(\bar u(\cd);v(\cd))=\dbE\int_0^T\Big\{\lan
b_x(T,s)^T\m(T)+\si_x(T,s)^T\n(s)+\int_0^Tf_x(r,s)^Tdr\\
\ns\ds\qq\qq+\(\psi_{x'}(s)^T
+\int_s^Tg_{x'}(s,r)^Tdr\)\xi(s)+\int_0^sg_x(r,s)^T\xi(r)
dr,X_1(s)\ran\\
\ns\ds\qq\qq+\lan b_u(T,s)^T\m(T)+\si_u(T,s)^T\n(s)+
\int_0^Tf_u(r,s)^Tdr+\int_0^sg_u(r,s)^T\xi(r)dr,v(s)-\bar
u(s)\ran\Big\}ds.\ea$$
Thus, $X_1(T)$ is eliminated.

\ms

\it Step 3. \rm Eliminate $X_1(\cd)$. We now apply Corollary 4.2
with
$$A_0(t,s)=b_x(t,s),\qq C_0(t,s)=\si_x(t,s),$$
and
$$\ba{ll}
\ns\ds\f(t)=\int_0^tb_u(t,s)[v(s)-\bar u(s)]ds+\int_0^t\si_u(t,s)[v(s)-\bar u(s)]dW(s),\\
\ns\ds\psi(t)=b_x(T,s)^T\m(T)+\si_x(T,s)^T\n(s)+\int_0^Tf_x(r,s)^Tdr\\
\ns\ds\qq\qq+\(\psi_{x'}(s)^T
+\int_s^Tg_{x'}(s,r)^Tdr\)\xi(s)+\int_0^sg_x(r,s)^T\xi(r)dr.\ea$$
Then we introduce the following BSVIE:
$$\ba{ll}
\ns\ds
p(t)=b_x(T,s)^T\m(T)+\si_x(T,s)^T\n(s)+\int_0^Tf_x(r,s)^Tdr\\
\ns\ds\qq\qq+\(\psi_{x'}(s)^T
+\int_s^Tg_{x'}(s,r)^Tdr\)\xi(s)+\int_0^sg_x(r,s)^T\xi(r)
dr\\
\ns\ds\qq\qq+\int_t^T\[b_x(s,t)^Tp(s)+\si_x(s,t)^Tq(s,t)\]ds-\int_t^Tq(t,s)
dW(s),\q t\in[0,T].\ea$$
Using Corollary 4.2, we have
$$\ba{ll}
\ns\ds\dbE\int_0^T\lan b_x(T,s)^T\m(T)+\si_x(T,s)^T\n(s)+\int_0^Tf_x(r,s)dr\\
\ns\ds\qq\qq+\(\psi_{x'}(s)^T
+\int_s^Tg_{x'}(s,r)^Tdr\)\xi(s)+\int_0^sg_x(r,s)^T\xi(r)dr,X_1(t)\ran
dt\\
\ns\ds\equiv\dbE\int_0^T\lan\psi(t),X_1(t)\ran
dt=\dbE\int_0^T\lan\f(t),p(t)\ran dt\\
\ns\ds=\dbE\int_0^T\[\int_0^t\lan b_u(t,s)[v(s)-\bar u(s)],p(t)\ran
ds+\int_0^t\lan\si_u(t,s)[v(s)-\bar u(s)],q(t,s)\ran ds\]dt.\ea$$
Hence,
$$\ba{ll}
\ns\ds J_2(\bar u(\cd);v(\cd))=\dbE\int_0^T\Big\{\lan
b_x(T,s)^T\m(T)+\si_x(T,s)^T\n(s)+\int_0^Tf_x(r,s)^Tdr\\
\ns\ds\qq\qq+\(\psi_{x'}(s)^T
+\int_s^Tg_{x'}(s,r)^Tdr\)\xi(s)+\int_0^sg_x(r,s)^T\xi(r)
dr,X_1(s)\ran\\
\ns\ds\qq\qq+\lan b_u(T,s)^T\m(T)+\si_u(T,s)^T\n(s)+
\int_0^Tf_u(r,s)^Tdr+\int_0^sg_u(r,s)^T\xi(r)dr,v(s)-\bar
u(s)\ran\Big\}ds\\
\ns\ds=\dbE\int_0^T\Big\{\lan b_u(T,s)^T\m(T)+\si_u(T,s)^T\n(s)+
\int_0^Tf_u(r,s)^Tdr\\
\ns\ds\qq\qq+\int_0^sg_u(r,s)^T\xi(r)dr+\int_s^Tb_u(r,s)^Tp(r)dr+
\int_s^T\si_u(r,s)^Tq(r,s)dr,v(s)-\bar u(s)\ran\Big\}ds.\ea$$
Then the maximum principle follows. \endpf

\ms

\section{Concluding Remarks}

In this paper, we have formulated two optimal control problems for
forward-backward stochastic Volterra integral equations.
Corresponding Pontryagin type maximum principles are established.
The key contribution of this paper is the discovery of a general
duality principle between a linear BSVIE and a linear
Fredholm-Volterra stochastic integral equation with mean-field. This
result substantially extends the corresponding result found in
\cite{Yong 2006,Yong 2008}. Such a duality principle enables us to
prove the maximum principles. It is known that when studying
stochastic integral equations, the It\^o's formula, which is very
powerful in studying SDEs, is absent. Therefore, one has to carry
out all relevant calculations without differentiation. We expect
that the techniques/ideas developed in this paper might have some
significant impacts on other type problems when It\^o's formula is
not directly applicable.


\begin{thebibliography}{00}
\rm

\bibitem{Ainslie 1992} G. Ainslie, Picoeconomics: The strategic
interaction of successive motivational states within the person.
Cambridge Univ. Press 1992.

\rm
%
%\bibitem{Cheridito-Delbaen-Kupper 2004} P. Cheridito, F. Delbaen, and M. Kupper, \it
%Coherent and covex monetary risk measures for bounded
%c$\grave{a}$dl$\grave{a}$g processes, \sl Stoch. process Appl. \rm
%112 (2004), 1-22.
%
%\rm
%\bibitem{Cheridito-Delbaen-Kupper 2006} P. Cheridito, F. Delbaen, and M. Kupper, \it
%Coherent and covex monetary risk measures for unbounded processes,
%\sl Finance Stoch. \rm 10 (2006), 427-448.
%
%\bibitem{Cheridito-Delbaen-Kupper 2006-2} P. Cheridito, F. Delbaen, and M. Kupper,
%\it Dynamic monetary risk measures for bounded discrete time
%processes, \sl Electron, J. Probab. \rm 11 (2006), 57-106.


\bibitem{Duffie-Epstein 1992} D. Duffie and L.G. Epstein, \it Stochastic differential utility,
\sl Econometrica, \rm 60 (1992), 353-394.

\bibitem{Duffie-Huang 1986} D. Duffie and C. F. Huang, \it Stochastic production-exchange
equilibria, \rm Reserach paper No.974, Graduate School of Business,
Stanford University, Stanford (1986).

\bibitem{Kamien-Muller 1976} M. I. Kamien and E. Muller, \it Optimal control with integral
state equations, \sl Rev. Eco. Study, \rm 43 (1976), 469--473.

\bibitem{Karatzas-Shreve 1988} I. Karatzas and S. E. Shreve, \it Brownian Motion and Stochastic Calculus.
Springer, Heidelberg (1988).

\rm
\bibitem{Lin 2002} J. Lin, \it Adapted solution of backward stochastic
nonlinear Volterra integral equation, \sl Stoch. Anal. Appl., \rm 20
(2002), 165--183.

\bibitem{Ma-Yong 1999} J. Ma and J. Yong, \it Forward-Backward Stochastic Differential Equations and Their Applications,
Springer-Verlag, Berlin. (1999) \rm

\bibitem{Pardoux-Protter 1990} E. Pardoux and P. Protter, \it Stochastic Volterra equations
with anticipating coefficients, \sl Ann. Probab., \rm 18 (1990),
1635--1655.

\bibitem{Peng 1990} S. Peng, \it A general stochastic maximum principle for
optimal control problems, \sl SIAM J. Control Optim., \rm 28 (1990),
966--979.

\bibitem{Peng 1993} S. Peng, \it Backward stochastic differential equations
and application to optimal control, \sl Appl. Math. Optim., \rm 27
(1993), 125--144.


\bibitem{Shi-Wang 2012} Y. Shi and T. Wang, \it Solvability of general backward
stochastic Volterra integral equations, \sl J. Korean. Math. Soc.,
\rm 49 (2012), 1301--1321.

\bibitem{Shi-Wang-Yong 2013} Y.~Shi, T.~Wang, and J.~Yong, \it
Mean-field backward stochastic Volterra integral equations, \sl Dis.
Cont. Dyn. System, Ser. B, \rm 18 (2013), 1929--1967.


%\bibitem{Vinokurov 1969} V. R. Vinokurov, \it Optimal control of processes described
%by integral equations, \sl SIAM J. Control, \rm 7 (1969), 324--355.

\bibitem{Wang-Shi 2010} T. Wang and Y. Shi, \it Symmetrical solutions of backward
stochastic Volterra integral equations and applications, \sl
Discrete Contin. Dyn. Syst. B., \rm 14 (2010), 251--274.

\bibitem{Wang-Yong 2013} T. Wang and J. Yong, \it Comparison theorems for backward
stochastic Volterra integral equations, \rm arxiv:1208.2064v1 [math.
PR] 10 Aug 2012.

\bibitem{Wang-Zhang 2007} Z. Wang and X. Zhang, \it Non-Lipschitz backward
stochastic volterra type equations with jumps, \sl Stoch. Dyn., \rm
7 (2007), 479--496.



\bibitem{Yong 2006} J. Yong, \it Backward stochastic Volterra integral
equations and some related problems, \sl Stoch. Proc. Appl., \rm 116
(2006), 779--795.

\bibitem{Yong 2007} J. Yong, \it Continuous-time dynamic risk measures by
backward stochastic Volterra integral equations, \sl Appl. Anal. \rm
86 (2007), 1429--1442.

\bibitem{Yong 2008} J. Yong, \it Well-posedness and regularity of
backward stochastic Volterra integral equation, \sl Probab. Theory
Relat. Fields, \rm 142 (2008), 21--77.

\bibitem{Yong-Zhou 1999} J. Yong and X. Y. Zhou, \sl Stochstic Control: Hamiltonian
Systems and HJB Equations, \rm Springer-Verlag, New York, 1999.



\end{thebibliography}
\end{document}